\documentclass[11pt]{article}
\usepackage{amsmath,amssymb}
\usepackage{graphics}
\usepackage{color}
\def\tto{\;{\lower 1pt \hbox{$\rightarrow$}}\kern -10pt
\hbox{\raise 2pt \hbox{$\rightarrow$}}\;}

\def\Tilde{\widetilde}

\def\lam{\lam}

\def\h{\hfill\Box}
\def\R{\mathbb{R}}
\def\N{\mathbb{N}}

\def\ox{\bar{x}}
\def\oy{\bar{y}}

\def\cone{\mbox{\rm cone}\,}

\def\dom{\mbox{\rm dom}\,}
\def\rge{\mbox{\rm rge}\,}

\def\conv{\mbox{\rm conv}\,}

\def\cl*co{\mbox{\rm cl}^*\mbox{\rm co}\,}
\def\cl{\mbox{\rm cl}\,}
\def\cl{\mbox{\rm cl}\,}
\def\vcl{\mbox{\rm vcl}\,}

\def\h{\hfill\triangle}

\def\O{\Omega}
\def\ph{\varphi}
\def\emp{\emptyset}
\def\st{\stackrel}

\def\ve{\varepsilon}
\def\Th{\Theta}

\def\D{\mathcal{D}\;\!}

\def\C{\mathcal{C}\;\!}

\definecolor{orange}{rgb}{1,0.5,0}

\begin{document}
\begin{center}
\vspace*{0.3in}{\bf\large VECTOR OPTIMIZATION WITH DOMINATION STRUCTURES:\\VARIATIONAL PRINCIPLES AND APPLICATIONS}
\\[2ex]
{\bf\large Truong Q. Bao\footnote{Department of Mathematics $\&$ Computer Science, Northern Michigan University, Marquette, Michigan 49855, USA (btruong@nmu.edu).}, Boris S. Mordukhovich\footnote{Department of Mathematics, Wayne State University, Detroit, Michigan 48202, USA (boris@math.wayne.edu). Research of this author was partly supported by the US National Science Foundation under grants DMS-1512846 and DMS-1808978, by the US Air Force Office of Scientific Research under grant \#15RT0462, and by the Australian Research Council Discovery Project DP-190100555.}, Antoine Soubeyran\footnote{Aix-Marseille University (Aix-Marseille School of Economics), CNRS \& EHESS, Marseille 13002, France (antoine.soubeyran@gmail.com).}\\
and Christiane Tammer\footnote{Martin-Luther-University Halle-Wittenberg, Faculty of Natural Sciences II, Institute of Mathematics, D-06099 Halle (Saale), Germany (christiane.tammer@mathematik.uni-halle.de).}}
\end{center}

\small {\bf Abstract}. This paper addresses a large class of vector optimization problems in infinite-dimensional spaces with respect to two important binary relations derived from domination structures. Motivated by theoretical challenges as well as by applications to some models in behavioral sciences, we establish new variational principles that can be viewed as far-going extensions of the Ekeland variational principle to cover domination vector settings. Our approach combines advantages of both primal and dual techniques in variational analysis with providing useful sufficient conditions for the existence of variational traps in behavioral science models with variable domination structures.\\[1ex]

\newtheorem{Theorem}{Theorem}[section]
\newtheorem{Proposition}[Theorem]{Proposition}
\newtheorem{Remark}[Theorem]{Remark}
\newtheorem{Lemma}[Theorem]{Lemma}
\newtheorem{Corollary}[Theorem]{Corollary}
\newtheorem{Definition}[Theorem]{Definition}
\newtheorem{Example}[Theorem]{Example}
\newtheorem{Comment}[Theorem]{Comment}
\newenvironment{Proof} {\noindent {\it Proof.\,  }\hspace{1pt}}
{\hspace{1pt}\hfill $\blacksquare$ \medskip }

\normalsize

\section{Introduction}\label{intro}

It has been well recognized over the years that problems of vector and set/set-valued optimization have great many mathematical challenges and important intrinsic issues, which significantly differ them from the conventional areas of scalar optimization. Thus such problems require developing novel tools of analysis to deal with their theory and applications. The spectrum of applications of vector and set-valued optimization is indeed enormous: economics, finance, ecology, radiotherapy treatment in medicine, environmental and behavioral sciences to name just a few; see \cite{bcs16,bms15,e14,grtz03,hamel,jahn,ktz,m18,ms19,s09} for more information and references. Variational principles, together with variational techniques and tools of generalized differentiation in set-valued and variational analysis, provide powerful machinery for the study and applications of vector and set optimality, particularly related to Pareto-type optimal/efficient solutions.  We refer the reader to, e.g., \cite{bm10,bms15b,grtz03,jahn,ktz,m06,m18,ms19} and the vast bibliographies therein for various approaches, concepts, and results in these and related directions.

This paper is devoted to the study and applications of vector optimization problems with {\em domination structures} that can be viewed as extensions of {\em variable} ordering cones. The {\em nondomination} concept for problems of multiobjective optimization (i.e., with finitely many scalar objectives) was introduced by Yu \cite{yu74} and then was studied and developed in many publications; see, e.g., \cite{bcy76,e14} and the references therein.  This concept is significantly more general than the conventional (Pareto) efficiency concept in vector optimization with fixed ordering cones; it has been realized as a crucial factor for a variety of applications to decision making, games, etc.

Here, we consider several solution concepts for general problems of vector optimization with domination structures while mainly focusing on two binary relations associated with a given domination/variable domination structure: {\em nondomination} and {\em efficiency}. After revealing important properties of both nondominated and efficient solutions with respect to variable domination structure in general linear space settings, we turn to the study of ordered-value mappings defined on {\em quasimetric} decision spaces. Besides mathematical novelty and interest, our motivation to involve quasimetric spaces (i.e., spaces with nonsymmetric distances) into consideration is due to unavoidable appearing such spaces in models of {\em behavioral sciences}; see the discussions in Section~\ref{behav}.\vspace*{0.03in}

The composition of the paper is as follows. In Section~\ref{sec:domin} we formulate vector optimization problems with preference relations  induced by the notion of {\em domination} defined via a set-valued mapping $\D:Y\rightrightarrows Y$ from the image space $Y$ to itself. This notion is a practically motivated extension of ordering structures given by {\em variable cones}. We define here several notions of {\em optimal solutions} with respect to two binary relations induced by a given domination structure, establish relationships between them, and then reveal some of their basic properties.

Section~\ref{overview} {\em overviews} and discusses known achievements in the theory of variational principles of the {\em Ekeland-type} for vector optimization problems with variable domination structures. For the reader's convenience and making comparison with our new developments obtained below, simplified proofs and clarifications of the major known results are given in this section.

Section~\ref{new} is the {\em culmination} of the paper, which presents new Ekeland-type {\em variational principles} in vector optimization with general domination structures and cost mappings defined on quasimetric spaces. We have two main motivations to develop these novel results. The {\em first motivation} comes from a strong mathematical call to obtain variational principles of the aforementioned type with taking into account drawbacks of the known results and their proofs discussed in Section~\ref{overview} as well as significant challenges that intrinsically appear in the new framework under consideration. The {\em second motivation} comes from the aimed applications to models of behavioral sciences in the vein of Soubeyran's {\em variational rationality} approach, which unavoidably involves quasimetric spaces (even in finite dimensions) and highly benefits from imposing domination structures on decision spaces; see \cite{s09,s19a,s19b} and the last section below. It is worth mentioning here that the progress achieved in this paper on variational principles and their applications is strongly based on the marriage of variational methods to Gerstewitz's {\em nonlinear scalarization} functional \cite{gers,ger-diss} and its recent developments given \cite{BouTam2017,gnrt16,w17}.

The concluding Section~\ref{behav} is devoted to {\em applications} of the obtained variational results in the general framework of vector optimization to some {\em behavioral science models} via developing the variational rationality approach to human dynamics initiated and conceptionally described in \cite{s09,s10,s16,s19a,s19b}. We first briefly review basic concepts of variational rationality that are closely related to quasimetric and domination structures. The major notions to analyze in these frameworks are {\em variational traps} of different types. Using our variational developments allows us to establish the existence of the so-called {\em ex ante} (before moving) and {\em ex post} (after moving) traps. In this way, we formulate
generalized efficiency and domination structures, which extend those in \cite{yu74} to the settings when resistance to move matters, and then derive the existence results in the ex ante and ex post traps in the new settings. Observe to this end that the proofs of the variational principles developed in Section~\ref{new}  provide {\em constructive dynamic procedures} to obtain such variational traps.\vspace*{-0.1in}

%%%%%%%%%%%%%%%%%%%%%%%%%%%%%%%%%%%%%%%%%%%%%%%%%%%%%%%%%%%%%%%%%%%%%%%%%%%%%%%

\section{Solution Concepts for Vector Optimization Problems with Respect to Domination Structures}\label{sec:domin}

This section is devoted to introducing the main concepts of our study,  establishing relationships between them, and revealing some of their important properties.

First we recall the classical notion of vector optimization with respect to a fixed ordering cone.

\begin{Definition}[\bf Pareto preorder]\label{def:Pareto} Let $Y$ be a linear space, let $C$ be a convex cone in $Y$, and let $y,v\in Y$. The {\sc Pareto preorder} on $Y$ denoted by $\leq_C$ is defined by
\begin{equation*}
v\leq_C y :\Longleftrightarrow v \in y - C \Longleftrightarrow y \in v + C.
\end{equation*}
\end{Definition}

Given two vectors $y$ and $v$ in a decision linear space $Y$, we write $v = y + d$ for some vector $d\in Y$. If $y$ is preferred by the decision maker to $v$, then $d$ can be viewed as a {\em domination factor}. The set of all the domination factors for $y$ together with the zero vector $0\in Y$ is denoted by $\D(y)$. Then, the multifunction $\D:Y\rightrightarrows Y$ is called a {\em domination structure}. It is also called a {\em variable ordering structure} in the majority of publications if $\D(y)$ is an {\em ordering cone} for each $y\in Y$.

The concept of domination structures was introduced by Yu in \cite{yu73dom}, where the sets $\D(y)$ are supposed to be cones. Yu defined a domination structure as a family of cones $\D(y)$, whereas Engau \cite{Engau08} considered it as a set-valued mapping. Domination factors were launched by Bergstresser, Charnes, and Yu \cite{bcy76} in a finite-dimensional setting with respect to convex domination sets.\vspace*{0.03in}

In contrast to vector optimization with a fixed ordering cone, we now define two binary relations in $Y$ with respect to the choice of domination sets. These binary relations will be used in the formulation of the variational principles in Sections~\ref{overview} and \ref{new} and also in the context of applications in behavioral sciences in Section~\ref{behav}.

\begin{Definition}[\bf binary relations]\label{def:binary} Given a domination structure $\D:Y\rightrightarrows Y$ in a linear space $Y$, and given vector $y,v\in Y$, we introduce the following binary relations:
\begin{itemize}
\item[\bf(i)] The {\sc nondomination binary relation} denoted by $\leq_{N,\D}$ is defined by
\begin{eqnarray*}
v \leq_{N,\D} y :\Longleftrightarrow y \in v + \D(v)\Longleftrightarrow B(y,v) := y - v \in \D(v).
\end{eqnarray*}
\item[\bf(ii)] The {\sc efficiency binary relation} denoted by $\leq_{E,\D}$ is defined by
\begin{eqnarray*}
v \leq_{E,\D} y :\Longleftrightarrow v \in y - \D(y)\Longleftrightarrow B(y,v) := y - v \in \D(y).
\end{eqnarray*}
\end{itemize}
\end{Definition}

The notation $B(y,v)$ allows us to unify the aforementioned binary relations being also useful in the application Section~\ref{behav}. It signifies the so-called ``worthwhile balance without inconvenience to move." The unified relation $B(y,v) \in \D(r)$, where the chosen reference point $r$ stands for $N$ as $r = v$, for $E$ as $r = y$, and for $\Theta$ as $r = v_0 = f(x_0)$. The worthwhile balance without inconvenience to move is $B(f(x),f(u)) = f(x) - f(u)$, the worthwhile balance without inconvenience to move and $r = N$ is $B(f(x)- \sqrt{\ve}q(x,u),f(u)) = f(x) - \sqrt{\ve}q(x,u) - f(u)$, and the worthwhile balance without inconvenience to move and $r = E$ is $B(f(x),f(u) - \sqrt{\ve}q(x,u)) = f(x) - f(u) - \sqrt{\ve}q(x,u)$. For simplicity, we drop the subscript $\D$ in the above binary notations if the context is clear. When $\D(y) \equiv C$ for some ordering cone of $Y$, both domination and efficiency binary relations reduce to the classical Pareto binary relation generated by $C$, i.e., $\leq_N = \leq_E = \leq_C$.\vspace*{0.05in}

Given further a mapping $f: X\to Y$ acting from a nonempty set to a linear space, we consider two solution concepts corresponding to both binary relations introduced in Definition~\ref{def:binary}. These notions are important to deriving of the variational principles in Sections~\ref{overview} and \ref{new} with the subsequent applications to the models of behavioral sciences given in Section~\ref{behav}. Denote in what follows $\dom f:=\{x\in X\;|\;f(x)\ne\emp\}$ and $\rge f := f(X)$ with $f(X):= \cup_{x \in X} f(x)$.

\begin{Definition}{\bf(nondominated and efficient solutions with respect to domination structures).}\label{def:sol} Let $f:X\to Y$ be a mapping from a nonempty set to a linear space, and let $\D:Y\rightrightarrows {Y}$ be a domination structure in the image space $Y$. Given $\ox\in\dom f$, we say that:
\begin{itemize}
\item[\bf(i)] $\ox$ is a {\sc conventional nondominated solution} of $f$ with respect to $\D$, or a conventional $\D$-nondominated solution, or a conventional $\leq_{N}$-minimal solution, if
$$
\forall x\in \dom f,\;f(x)\leq_N f(\ox) \Longrightarrow f(\ox)\leq_N f(x).
$$
\item[\bf(ii)] $\ox$ is a {\sc $\D$-nondominated solution} of $f$ with respect to $\D$ if
$$
\forall x\in \dom f,\;f(x) \neq f(\ox) \Longrightarrow f(x)\not\leq_N f(\ox).
$$
\item[\bf(iii)] $\ox$ is a {\sc conventional efficient solution} of $f$ with respect to $\D$, or a conventional $\D$-efficient solution, or a conventional $\leq_{E}$-minimal solution, if
$$
\forall x\in \dom f,\;f(x)\leq_E f(\ox) \Longrightarrow f(\ox)\leq_E f(x).
$$
\item[\bf(iv)] $\ox$ is a {\sc $\D$-efficient solution} of $f$ with respect to $\D$ if
$$
\forall x\in \dom f,\;f(x) \neq f(\ox) \Longrightarrow f(x)\not\leq_E f(\ox),
$$
which is equivalent to the condition
$$
\rge{f}\cap(f(\ox) - \D(f(\ox))) = \{f(\ox)\}.
$$
\end{itemize}
\end{Definition}

Recall that the concept of $\D$-nondominated solutions in Definition \ref{def:sol}(ii) was initiated by Yu \cite{yu73dom,yu74} for conic domination structures. The concept of $\D$-efficient solutions of $f$ with respect to $\D$ in Definition \ref{def:sol}(iv) was introduced by Chen, Huang, and Yang \cite[Definition 1.13]{CheHuaYan05} under the name of ``nondominated-like minimal points"; see also Chen and Yang \cite[Definition 3.1]{CheYan02}. We will deal with (approximate) $\D$-efficient and $\D$-nondominated solutions in Sections~\ref{overview}, \ref{new}, and \ref{behav}.

Observe also that a $\D$-efficient solution $\ox$ of $f$ is an element, which is not dominated by another point $x$ with respect to the associated set $\D (f(\overline{x}))$ at the $\D$-efficient solution $\overline{x}$. However, given a $\D$-nondominated solution $\overline{x}$ of $f$, a domination set $\D (f(x))$ is a set associated with another point $x$. Important properties of these elements can be found in \cite{CheYan02, CheHuaYan05,e14,Eic13,Eic12,yu74}.\vspace*{0.05in}

In order to combine these two solution concepts, we use in the following proposition the language of {\em $\leq_\circ$-minimality}, where $\leq_\circ$ stands for either the domination binary relation $\leq_N$, or for the efficient binary relation $\leq_E$ formulated in Definition~\ref{def:binary}.

\begin{Proposition}[\bf relationships between minimal solutions, I]\label{relat1} Let `$\leq_\circ$' stand for both the domination binary relation $\leq_N$ and the efficient binary relation $\leq_E$ taken from Definition~{\rm\ref{def:binary}}. Then, we have the relationships:
\begin{itemize}
\item[\bf(i)] If $\ox$ is a $\leq_\circ$-minimal solution of $f$, then it is a conventional $\leq_\circ$-minimal solution of $f$.

\item[\bf(ii)] Assume that the pointedness condition for $\{f,\D\}$ at $\ox\in\dom f$
\begin{eqnarray}\label{cond:pointed}
\forall x\in \dom{f},\;\D(f(x))\cap(-\D(f(\ox))) = \{0\}
\end{eqnarray}
holds. If $\ox$ is a conventional $\leq_\circ$-minimal solution of $f$, then it is $\leq_\circ$-minimal to $f$.
\end{itemize}
\end{Proposition}
{\bf Proof}. Let us verify both conclusions in this proposition for the case of nondominated solutions; the proof for efficient solutions is similar.

To justify (i), assume that $\ox$ is a $\D$-nondominated solution of $f$, i.e.,
$$
\forall x\in \dom f,\;f(x) \neq f(\ox) \Longrightarrow f(x)\not\leq_N f(\ox),
$$
which is equivalent to the implication
$$
\forall x\in \dom{f},\;f(x)\leq_N f(\ox)  \Longrightarrow f(x) = f(\ox).
$$
By $f(x) = f(\ox)$, we have $f(\ox)\leq_N f(x)$, and so $\ox$ is a conventional $\D$-nondominated solution of $f$.

To prove (ii), assume that (\ref{cond:pointed}) is satisfied, and that $\ox$ is a conventional $\D$-nondominated solution of $f$. To check the $\D$-nondomination of $\ox$ to $f$, fix an arbitrary element $x\in \dom{f}$ satisfying $f(x) \neq f(\ox)$. We claim that $f(x) \not\leq_N f(\ox)$. Arguing by contraposition, suppose that  $f(x) \leq_N f(\ox)$. The conventional $\D$-nondominatedness of $\ox$ to $f$ yields $f(\ox)\leq_N f(x)$. Then, Definition~\ref{def:binary}(i) tells us that
$$
f(\ox) \in f(x) + \D(f(x))\;\mbox{\rm and } f(x) \in f(\ox) + \D(f(\ox)),
$$
and therefore $f(\ox) - f(x) \in \D(f(x)) \cap(- \D(f(\ox))) = \{0\}$, where the last equality holds due to the pointedness condition for $\{f,\D\}$ at $\ox$. Thus we get  $f(\ox)=f(x)$, a contradiction, which shows that $f(x) \not\leq_N f(\ox)$. Since $x$ was chosen arbitrary in $\dom(f)$ while satisfying $f(x) \neq f(\ox)$, we verify that $\ox$ is a $\D$-nondominated solution of $f$ and hence complete the proof of the proposition. $\h$\vspace*{0.05in}

Observe that when $\D(y) \equiv \Theta$ is a fixed domination set or $\D(y) \equiv C$ is a fixed ordering cone, there is no difference between the two concepts of nondomination and efficiency defined in Definition~\ref{def:sol}, and they both reduce to {\em Pareto minimality}. In such a situation, the pointedness condition for $\D$ and $f$ at $\ox\in\dom f$ is nothing but the pointedness property of the ordering set $\Th$ and of the ordering cone $C$, respectively.\vspace*{0.05in}

Next we establish relationships between minimal solutions with respect to a {\em domination structure} $\D$ and Pareto minimal solutions with respect to a {\em fixed domination set} $\Th$. Denote the two domination sets associated with $\{f,\D\}$ by
$$
 {\Theta^{u}_{\D}} := \bigcup\big\{ \D(f(x))\;\big|\;x\in \dom{f}\big\} \quad \mbox{\rm and } \quad  {\Theta^{i}_{\D}} := \bigcap\big\{ \D(f(x))\;\big|\;x\in  {\dom f }\big\}
$$
and call them the {\em union} (respectively, {\em intersection}) {\em domination set} for $f$ and $\D$. We skip mentioning $f$ in the above notations for simplicity.

\begin{Proposition}[\bf relationships between minimal solutions, II]\label{relst2} The following hold:
\begin{itemize}
\item[\bf(i)] If $\ox$ is a conventional $\D$-efficient solution of $f$, then it is $\Tilde{\D}$-efficient to $f$, where $\Tilde{\D}:Y\rightrightarrows {Y}$ is defined by $\Tilde{\D}(y) := \D(\oy)\setminus(-\Theta^u_{\D})$.

\item[\bf(ii)] If $\ox$ is a $\D$-efficient solution of $f$, then it is $\Th$-minimal with $\Th = \D(f(\ox))$.

\item[\bf(iii)] If $\ox$ is a $\D$-nondominated solution of $f$, then it is $ {\Theta^{i}_{\D}}$-minimal to $f$.

\item[\bf(iv)] If $\ox$ is a $ {\Theta^{u}_{\D}}$-minimal solution of $f$, then it is $\D$-nondominated to $f$.
\end{itemize}
\end{Proposition}
{\bf Proof}. (i) Assume that $\ox$ is a conventional $\D$-efficient solution of $f$, i.e.,
\begin{eqnarray*}
& & \forall x\in \dom{f},\;f(x)\leq_{E,\D} f(\ox) \Longrightarrow f(\ox)\leq_{E,\D} f(x)\\
& \Longleftrightarrow & \forall x\in \dom{f},\;f(x)\in f(\ox) - \D(f(\ox)) \Longrightarrow f(\ox) \in  f(x) - \D(f(\ox)).
\end{eqnarray*}
Arguing by contraposition, suppose that $\ox$ is not a $\Tilde{\D}$-efficient solution of $f$. Then, we could find $x\in \dom{f}$ satisfying $f(x) \neq f(\ox)$ and such that
$f(x) \leq_{E,\Tilde{\D}} f(\ox)$, i.e.,
\begin{eqnarray}\label{cond:tD}
f(x) \in f(\ox) - \Tilde{\D}(f(\ox)) = f(\ox) - \D(f(\ox))\setminus(-\Theta^u_{\D}) \subseteq f(\ox) - \D(f(\ox)),
\end{eqnarray}
which clearly implies that $f(x) - f(\ox) \not\in \Theta^u_{\D}$ and $f(x) \leq_{E,\D} f(\ox)$. Since $\ox$ is a conventional $\D$-efficient solution of $f$, we have $f(\ox) \leq_{E,\D} f(x)$, i.e.,
$$
f(\ox) \in f(x) - \D(f(x)) \Longrightarrow f(x) - f(\ox) \in \D(f(x)) \subseteq  {\Theta^{u}_{\D}}.
$$
The obtained contradiction verifies the implication in (i).

(ii) This is straightforward from the definitions. Indeed, we have
\begin{eqnarray*}
& & \mbox{$\ox$ is a $\D$-efficient solution of $f$}\\
&\Longleftrightarrow & \forall x\in \dom{f},\;f(x) \neq f(\ox) \Longrightarrow f(x)\not\leq_E f(\ox)\\
&\Longleftrightarrow & \forall x\in \dom{f},\;f(x) \neq f(\ox) \Longrightarrow f(x)\not\in f(\ox) - \D(f(\ox))\\
&\Longleftrightarrow & \forall x\in \dom{f},\;f(x) \neq f(\ox) \Longrightarrow f(x)\not\leq_{\D(f(\ox))}   f(\ox)\\
&\Longleftrightarrow & \mbox{$\ox$ is a $\leq_{\D(f(\ox))}$-minimal solution of $f$}.
\end{eqnarray*}

(iii) This also follows from the definitions. Indeed, we have
\begin{eqnarray*}
& & \mbox{$\ox$ is a $\D$-nondominated solution of $f$}\\
&\Longleftrightarrow & \forall x\in \dom{f},\;f(x) \neq f(\ox) \Longrightarrow f(x)\not\leq_N f(\ox)\\
&\Longleftrightarrow & \forall x\in \dom{f},\;f(x) \neq f(\ox) \Longrightarrow f(\ox)\not\in f(x) + \D(f(x))\\
&\Longleftrightarrow & \forall x\in \dom{f},\;f(x) \neq f(\ox) \Longrightarrow f(\ox)- f(x) \in Y\setminus\D(f(x)) \subseteq Y\setminus {\Theta^{i}_{\D}}\\
&\Longleftarrow & \forall x\in \dom{f},\;f(x) \neq f(\ox) \Longrightarrow f(x)\not\leq_{ {\Theta^{i}_{\D}}}  f(\ox)\\
&\Longleftrightarrow & \mbox{$\ox$ is a (Pareto) $\leq_{ {\Theta^{i}_{\D}}}$-minimal solution of $f$}.
\end{eqnarray*}

(iv) Similarly to the above we get the equivalences
\begin{eqnarray*}
& & \mbox{$\ox$ is a $ {\Theta^{u}_{\D}}$-minimal solution of $f$}\\
&\Longleftrightarrow & \forall x\in \dom{f},\;f(x) \neq f(\ox) \Longrightarrow f(x)\not\leq_{ {\Theta^{u}_{\D}}}  f(\ox)\\
&\Longleftrightarrow & \forall x\in \dom{f},\;f(x) \neq f(\ox) \Longrightarrow f(\ox)\not\in f(x) +  {\Theta^{u}_{\D}}\\
&\Longleftrightarrow & \forall x\in \dom{f},\;f(x) \neq f(\ox) \Longrightarrow f(\ox)- f(x) \in Y\setminus {\Theta^{u}_{\D}} \subseteq Y\setminus\D(f(x))\\
&\Longleftarrow & \forall x\in \dom{f},\;f(x) \neq f(\ox) \Longrightarrow f(x) \not\leq_{N} f(\ox)\\
&\Longleftrightarrow & \mbox{$\ox$ is a $\D$-nondominated solution of $f$},
\end{eqnarray*}
which therefore complete the proof of the proposition. $\h$\vspace*{0.07in}

We illustrate the differences of solution notions above by the following example.

\begin{Example}[\bf differences between solution notions]\label{exa1} {\rm Let $X := \{1,2,3\}$, and let the values of a mapping $f:X\to \R^2$ are given by
$$
f(1) := A_1 = (0,0),\;f(2) := A_2 = (4,-2),\; \mbox{ and }\; f(3) := A_3 = (-2,1).
$$
Set $e_1 := (1,0)$, $e_2 := (0,1)$, $e_3 := (1,1)$, and $e_4 := (-1,-1)$ and then consider a domination structure $\D$ on $f(X)$ with the following values in $\R^2$:
$$
\D(A_1) := \conv\cone\{e_1,e_2\},  \D(A_2) := \conv\cone\{e_1,e_3,e_4\}, \; \mbox{ and }\; \D(A_3) := \conv\cone\{-e_2,e_3\}.
$$
Then, we have the optimal solutions:
\begin{itemize}
\item[a)] $1$ and $3$ are $\D$-efficient solutions of $f$.
\item[b)] $3$ is a $\D$-nondominated solution of $f$.
\item[c)] $3$ is a $\Theta^{i}_{\D}$-minimal solution of $f$.
\end{itemize}}
\end{Example}\vspace*{-0.15in}

%%%%%%%%%%%%%%%%%%%%%%%%%%%%%%%%%%%%%%%%%%%%%%%%%%%%%%%%%%%%%%%%%%%%%%%%%%%%%%%%%%%%%%%

\section{Overview and Elaborations of Known Results}\label{overview}

It has been well recognized that the {\em Ekeland variational principle} (EVP) plays a fundamental role in variational analysis and in a vast variety of applications, including those to vector and set-valued optimization; see, e.g., the books \cite{ktz,m06,m18} with the references and commentaries therein. Quite recently, several extensions of the EVP have been developed for problems of vector optimization with {\em domination structures}. Let us recall and elaborate them in this section, which makes a bridge to our new developments and applications in the subsequent sections.

In the first part of this section we deal with approximate solutions of $f$ as an extension of $\D$-efficient solutions of $f$ in the sense of Definition~\ref{def:sol}(iv). This concept of (weakly) approximate efficient solutions with respect to $\D$ is given in the next definition.

\begin{Definition}[\bf approximate efficient solutions]\label{d-appr-eff}
Let $f : X \to Y$, $\D : Y \rightrightarrows  Y$, $k\in Y\setminus\{0\}$, and $\ve\geq{0}$. Then, we have:
\begin{itemize}
\item[\bf(a)] An element $x_\ve \in X $ is called an {\sc $\ve{k}$-efficient solution} of $f$ with respect to $\D$ if
$$
{f}(X) \cap  (f(x_\ve) - \ve{k} - (\D(f(x_\ve))\setminus \{0\})) = \emptyset.
$$
\item[\bf(b)] An element $x_\ve \in X $ is called a {\sc weakly $\ve{k}$-efficient solution} of $f$ with respect to $\D$ if $ \operatorname*{int} \D(f(x_\ve)) \neq \emptyset$ and
$$
{f}(X) \cap  (f(x_\ve) - \ve{k} - \operatorname*{int} \D(f(x_\ve))) = \emptyset.
$$
\end{itemize}
\end{Definition}
Note that for the special case where $\ve = {0}$, the notion of (weakly) $\ve{k}$-efficient solutions of $f$ with respect to $\D$ reduces to that of $\D$-efficient solutions of $f$ formulated to Definition~\ref{def:sol}(iv).\vspace*{0.03in}

The classical EVP concerns approximate solutions of scalar optimization problems with extended-real valued, lower semicontinuous, and bounded from below objectives in the setting of complete metric spaces. First we recall extensions of EVP to the case of $\ve{k}$-efficient solutions of $f$ with respect to a domination structure $\D: Y \rightrightarrows  Y$ in the sense of Definition \ref{d-appr-eff}.\vspace*{0.05in}

In \cite[Theorem~3.1]{bms15}, Bao et al. established a version of EVP for set-valued mappings with ordering structures in quasimetric spaces. They used a variational approach based on an extended version of the Dancs-Hegedu\c{s}-Medvegyev\c{s} fixed point theorem. A simplified version of this result for vector-valued mappings from a complete metric space to a normed space is given below.

\begin{Theorem}[\bf EVP for conic domination structures]\label{thm:bms15}
Let $(X,d)$ be a complete metric space, let $Y$ be a normed space, and let $k\in Y\setminus\{0\}$. Given $f : X \to Y$ and a cone $\Theta\subset Y$. Consider a domination structure $\D : Y \rightrightarrows  Y$ such that the sets $\D(y)$ are proper, pointed, and closed cones for all $y\in\rge{f}$. Consider $\Theta_{\D}^i = \cap\{\D(f(x))\;|\; {x\in \dom{f}}\}$. Picking $\ve >0$, take an $\ve{k}$-efficient solution $x_\ve \in X$  of $f$ with respect to $\D$ together with $k\in\Theta_{\D}^i\setminus(-\Theta - \D(f(x_0)))$ and assume that:
\begin{itemize}
\item[\bf(A1)] {\sc $($boundedness condition$)$} $f$ is quasibounded with respect to $\Theta$ in the sense that there is a bounded subset of $Y$ such that $f(x) \subseteq M + \Theta$ for all $x\in \dom{f}$.

\item[\bf(A2)] {\sc $($limiting monotonicity condition$)$} $f$ satisfies the limiting decreasing continuity condition over $X$ with respect to $\D$ in the sense that for every sequence $\{x_n\}\subseteq X$ such that
$$
x_n \to x_\ast \mbox{ \; and\; }    {f(x_n) - f(x_{n+1}) \in \D(f(x_n))}
$$
we have  {$f(x_n) - f(x_{\ast}) \in \D(f(x_n))$} for all $n\in \mathbb{N}$.

\item[\bf(A3)] {\sc (transitivity condition for $\leq_E$)} $\D$ enjoys the monotonicity property on $\rge{f}$ in the sense of the application
$$
\forall x,u\in X,\;{f(u) - f(x) \in \D( f(u))} \Longrightarrow D(f(x)) \subseteq \D(f(u)).
$$
\end{itemize}
Then, there exists $\ox \in \dom{f}$ such that the following conditions hold:
\begin{itemize}
\item[\bf(i)] $\ox$ is an ${\ve}k$-efficient solution of $f$ with respect to $\D$.
\item[\bf(ii)] $d(x_\ve, \ox) \leq \sqrt{\ve}$.
\item[\bf(iii)] $\ox$ is an efficient solution of $f_{\ox}$ with respect to $\D$, where $f_{\ox} := f + \sqrt{\ve}d(\ox,\cdot)k$, i.e.,
$$\forall x\in X\setminus\{\ox\},   {f(\ox) - f(x) - \sqrt{\ve}d(\ox,x)k \not\in \D(f(\ox)).}$$
\end{itemize}
\end{Theorem}

\begin{Comment} {\bf(a)} {\rm Conclusions (i) and (ii) can be formulated in the form
$$  {f(x_\ve) - f(\ox) - \sqrt{\ve}d(x_\ve,\ox)k \in \D(f(x_\ve)),}$$
which clearly implies that both (i) and (ii) hold.

{\bf(b)} The requirement $k\in\Theta_{\D}^i$ is equivalent to that $k\in \D(f(x))$ for all $x\in \dom{f}$, and the monotonicity property of the variable ordering structure $\D$ is essential for the transitivity property of the efficiency binary relation $\leq_{E}$. Furthermore, both these assumptions are essential in the variational approach in \cite{bm10}.

{\bf(c)} The space $Y$ should be a normed space in \cite[Theorem~3.1]{bms15} instead of a real topological space, since the existence of a bounded set in condition (A1) is not defined in the latter case}.
\end{Comment}

In \cite[Theorem~5.1]{s14}, Soleimani established a version of EVP in vector optimization with a domination structure $\D: X \rightrightarrows  Y$ whose domination sets are {\em not necessarily cones}. He used the {\em scalarization approach} first developed in \cite{t92} for vector optimization with domination sets, i.e., in the case of a constant domination structure. The result formulated in Banach spaces holds in this setting. It was improved in \cite[Theorem~3.8]{best17}, while in \cite[Theorem~5.1]{s14} and \cite[Theorem~3.8]{best17} it was supposed that $X$ is a Banach space and that $Y$ a topological linear space.\vspace*{0.05in}

The Ekeland-type variational principle given in the next theorem is derived for $\ve{k}$-efficient solutions of $f$ with respect to $\D : Y \rightrightarrows  Y$ in the sense of Definition~\ref{d-appr-eff} under the assumption that $(X,d)$ is a complete metric space and that $Y$ a topological linear space. It provides a certain improvement of \cite[Theorem~5.1]{s14} and \cite[Theorem~3.8]{best17} with a simplified proof as presented below.

\begin{Theorem}[\bf EVP with relaxed conic domination]\label{thm:benahn} Let $(X,d)$ be a complete metric space, let $Y$ be a topological linear space, and let $k\in Y\setminus\{0\}$. Given a domination structure $\D : Y \rightrightarrows  Y$, a mapping $f : X \to Y$, and a number $\ve\geq{0}$, we consider an $\ve{k}$-efficient solution $x_\ve$ of $f$ with respect to $\D$ and denote $y_\ve:= f(x_\ve)$.
Assume the following conditions:
\begin{itemize}
\item[\bf(B1)] {\sc(boundedness condition)} $f$ is bounded from below in the sense that there is $\underline{y}\in Y$ such that $\forall x\in X , f(x) \in \underline{y} - \D(y_\ve)$.

\item[\bf(B2)] {\sc(lower semicontinuity condition)} $f$ is $(k,\D)$-lsc in the sense that for every $y\in f(X) $ and for every $t\in \mathbb{R}$ the sets
$$
M(y,t) :=\big\{x \in X\;\big|\;f(x)\in tk - \cl(\D(y))\big\}
$$
are closed in $X$.

\item[\bf(B3)] {\sc(scalarization conditions)}
\begin{itemize}
\item [\bf(B3-a)] $0 \in \D(y_\ve)$, $\D(y_\ve)$ is a proper, pointed, closed, and solid set with $\D(y_\ve) + \D(y_\ve) \subseteq \D(y_\ve)$ and $\D(y_\ve)+(0,\infty)k \subseteq \operatorname*{int} (\D(y_\ve))$.
\item [\bf(B3-b)] There is a cone-valued mapping $\C : Y \rightrightarrows Y$ satisfying $k\in \operatorname*{int}(\C(y_\ve))$ and $\D(y_\ve)+ (\C(y_\ve)\setminus\{0\}) \subseteq \operatorname*{int} \D(y_\ve)$.
\item [\bf(B3-c)] $\D(y) \subseteq \D(y_\ve)$ for all $y \leq_{\D(y_\ve)} y_\ve$.
\end{itemize}
\end{itemize}
Then, there exists $\ox\in \dom{f}$  such that we have the assertions:
\begin{itemize}
\item[\bf(i)] $\ox$  is an $\ve{k}$-efficient solution of $f$ with respect to $\C$, i.e.,
$$
\forall x\in X, \;    {f(\ox) - f(x) -  \sqrt{{\ve}}k \not\in  \C(\oy) \setminus \{0\},}
$$
where $\oy := f(\ox )$.

\item[\bf(ii)] $d(x_\ve,\ox) \leq \sqrt{\ve}$.

\item[\bf(iii)] $\ox$ is an efficient solution of the perturbed function $f_{\ox}$ with respect to $\C$, where $f_{\ox} := f + \sqrt{\ve}d(\ox,\cdot)k$, i.e.,
$$
\forall x\in  X,\;f(\ox) - f(x) - \sqrt{\ve}d(\ox,x)k \not\in  \C(\oy) \setminus \{0\}.
$$
\end{itemize}
\end{Theorem}
{\bf Proof}. Set  {$\Theta := \D(f(x_\ve))$}, define
$$
S := \big\{x\in  X\;\big|\;   {f(x_\ve) - f(x) - \sqrt{\ve}d( x_\ve,x)k \in \D(f(x_\ve))}\big\},
$$
and consider the scalarization function $\ph: Y\to \mathbb{R}$ given by
$$
\ph(y) := \ph_{\Th,k}(y) =\big\{t\in \mathbb{R}\;\big|\; y \in tk - \Theta\big\}.
$$
By \cite[Theorem~2.3.1]{grtz03}, condition (B3-a) ensures that $\dom{\ph} = Y$ and that $\ph$ is translation invariant along $k$. Furthermore, condition (B3-b) implies that $\ph$ is strictly $C(y_\ve)$-monotone in the sense that
$$
a \in b -(C(y_\ve) \setminus \{0\})\mbox{\;and\;} a\neq{b} \Longrightarrow \ph(a) < \ph(b).
$$
It is easy to check that the $\ve{k}$-efficiency to $f$ of $x_\ve$ guarantees that $x_\ve$ is an $\ve$-minimal solution of the scalarized function $\psi := \ph\circ{f}$.

By (B1) and (B2), the scalarized function $\psi$ is bounded from below and lower semicontinuous. The classical EVP ensures the existence of $\ox \in \dom{f}$ such that
\begin{itemize}
\item[\rm(i')] $\psi(\ox) + \sqrt{\ve}d(x_\ve,\ox) \leq \psi(x_\ve)$ and
\item[\rm(ii')] $\forall x\in  X \setminus\{\ox\},\;\psi(x) + \sqrt{\ve}d(\ox,x) > \psi(\ox)$.
\end{itemize}
Since $x_\ve$ is an $\ve$-minimal solution of $\psi$,  it follows from (i') that $d(x_\ve,\ox) \leq \sqrt{\ve}$. By \cite[Theorem~2.3.1]{grtz03} we have that $f(\ox) +\sqrt{\ve}d(x_\ve,\ox) k\in  y_\ve - \D(y_\ve)$, and thus $f(\ox)  \in  y_\ve - \operatorname*{int}\D(y_\ve)$ due to (B3-a). Arguing by contraposition, now we verify the fulfillment of (i). Indeed, suppose that $\ox$  is not an $\ve{k}$-efficient solution of $f$ with respect to $\C$, i.e., there is some $x\in X$ such that
$$f(x) +  \ve{k} \in f(\ox) - ( \C(\oy) \setminus \{0\}).
$$
Combining the two previous inclusions tells us that
\begin{eqnarray*}
f(x) +\ve{k} & \in & f(x_\ve) - \operatorname*{int}\D(y_\ve) - ( \C(\oy) \setminus \{0\}) \\
& \st{(B3-b)}{\subseteq} &  f(x_\ve) - \operatorname*{int}\D(y_\ve) - ( \D(\oy) \setminus \{0\})\\
& \st{(B3-c)}{\subseteq} & f(x_\ve) - \operatorname*{int}\D(y_\ve) - ( \D(y_\ve) \setminus \{0\})\\
& \subseteq & f(x_\ve) - ( \D(y_\ve) \setminus \{0\}),
\end{eqnarray*}
which clearly contradicts the $\ve{k}$-efficiency of $x_\ve$ for $f$ with respect to $\D$.

We complete the proof by verifying that (ii') yields (iii). Indeed, it follows that
\begin{eqnarray*}
(ii')& \Longleftrightarrow &  \forall x\in X\setminus\{\ox\},\;\ph(f(x) + \sqrt{\ve}d(\ox,x)k > \ph(f(\ox))\\
& \Longrightarrow & \forall x\in X\setminus\{\ox\},\;f(x) + \sqrt{\ve}d(\ox,x)k \not\in f(\ox)) -\D(f(x_\ve )) \\
& \Longrightarrow & \forall x\in X\setminus\{\ox\},\;f(x) + \sqrt{\ve}d(\ox,x)k \not\in f(\ox)) - (\C(f (\ox )) \setminus \{0\}),
\end{eqnarray*}
where the first implication holds due to \cite[Theorem~2.3.1]{grtz03} and the last one holds due to (B3-c). Thus (iii) is satisfied, which ends the proof of the theorem. $\h$

\begin{Comment} {\bf(a)} {\rm The result of \cite[Theorem~5.1]{s14} was formulated for ordering structures $\D_{X} : X\rightrightarrows  Y$ acting from the domain space to the image space of the mapping $f$. By using the same line as in the proof of Theorem~\ref{thm:benahn}, a better result can be established since the hypotheses of \cite[Theorem~5.1]{s14} are more restrictive. In particular, (B3) is assumed therein for all $x \in \dom{f}$ instead of at the element $x_\ve$. The idea of the proof in \cite[Theorem~5.1]{s14} is to scalarize the vector-valued mapping by using the nonlinear scalarized function $\ph := \ph_{\D(y_\ve),k}$. Thus it is sufficient to impose the assumptions on the domination set $\D(y_\ve)$ as in our elaborated proof.

{\bf(b)} The lower semicontinuity condition (B2) can be weaken to strictly decreasing lower semicontinuity as in \cite{bcs16}.

{\bf(c)} The boundedness condition of the scalar function $\psi$ is equivalent to the existence of a real number $m$ such that
\begin{eqnarray}\label{new-bounded}
\psi(x) = \ph(f(x)) > m \; \mbox{ for all } \; x\in \dom{f}\quad
\Longleftrightarrow \quad f(x) \not\in m k - \D(y_\ve).
\end{eqnarray}
It is easy to check that this condition is weaker than (B1).

{\bf(d)} Condition (B3) does not guarantee that the binary relation $\leq_{\D(y_\ve)}$ is transitive.

{\bf(e)} The result was formulated with respect to the cone-valued domination (i.e., ordering) structure $\C$. In the next section, we will establish a corresponding (even better) one in terms of a given domination structure $\D$, which may be nonconic.}
\end{Comment}\vspace*{-0.1in}

Next we recall Theorem~3.12 in \cite{best17}, which is a version of EVP with nonsolid domination sets for $\ve{k}$-efficient solutions of $f$ with respect to $\D: Y \rightrightarrows  Y$ in the sense of Definition \ref{d-appr-eff}. It was established by using the nonlinear scalarization approach.

\begin{Theorem}[\bf EVP with nonsolid domination sets]\label{t-nonsolid} Let $(X,d)$ be a complete metric space, let $Y$ be a Banach space, and let $k\in Y\setminus\{0\}$. Given a domination structure $\D : Y \rightrightarrows  Y$ and a mapping $f : X \to Y$, for each $\ve\geq{0}$ and $k\in Y\setminus\{0\}$ consider an $\ve{k}$-efficient solution $x_\ve$ of $f$ with respect to $\D$ and denote $y_\ve:= f(x_\ve)$. Impose the following assumptions:
\begin{itemize}
\item[\bf(C1)] {\sc (quasiboundedness condition)} $f$ is quasibounded from below with respect to $\D(y_\ve)$ in the sense that there is a bounded set $M \subseteq Y$ such that $f(x) \in M - \D(y_\ve)$ for all $x\in \dom{f}$.

\item[\bf(C2)] {\sc (lower continuity condition)} $f$ is $\D(y_\ve)$-lower semicontinuous over $\dom{f}$ in the sense that the level sets
$$
{\rm lev}(y;f) :=\big\{x\in \dom{f}\;\big|\;f(x) \in y - \D(y_\ve)\big\}
$$
are closed in $X$ for all $y\in Y$.

\item[\bf(C3)] {\sc (scalarization conditions)}
\begin{itemize}
\item[\bf(C1-a)] $\D(y_\ve)$ is a proper, closed, convex, and pointed cone.
\item[\bf(C2-b)] $\D(f(x)) \subseteq \D(y_\ve)$ for all $x\in \dom{f}$ with $ d(x_\ve,x) \le \sqrt{\ve}$.
\end{itemize}
\end{itemize}
Then, there exists an element $\ox \in \operatorname*{dom}f$ for which we have:
\begin{itemize}
\item[\bf(i)]  {$f(x_\ve) - f(\ox) \in  \D(y_\ve)$}, and thus $\ox$ is an $\ve{k}$-efficient solution of $f$ with respect to $\D$.
\item[\bf(ii)] $d(x_\ve,\ox) \leq \sqrt{\ve}$.
\item[\bf(iii)] $\ox$ is an efficient solution of $f_{\ox}$ with respect to $\D$, where  $f_{\ox} := f + \sqrt{\ve}d(\ox,\cdot)k$, i.e.,
$$
{f(\ox) - f(x) - \sqrt{\ve}d(\ox,x)k \not\in \D(f(\ox))}  \; \mbox{ for all } \; x\in \dom{f}.
$$
\end{itemize}
\end{Theorem}

\begin{Comment} {\bf(a)} {\rm Theorem~3.12 in \cite{best17} was formulated for $\ve{k}$-efficient solutions of constrained problems, where the cost mapping $f$ acted between two Banach spaces. However, the proof of the theorem holds true when the domain space is a complete metric one.

{\bf(b)} The quasiboundedness condition (C1) is more restrictive then condition \eqref{new-bounded}, which is equivalent to the boundedness from below of the scalarized function of $f$.

{\bf(c)} Condition (C2) implies that the composite function $\psi = \ph\circ{f}$ defined in the proof of Theorem~\ref{thm:benahn} is lower semicontinuous. Therefore, we can weaken it to the requirement that $\psi$ is strictly decreasingly lower semicontinuous.}
\end{Comment}

In contrast to the aforementioned developments of EVP for $\ve{k}$-efficient solutions of vector-valued mappings with respect to domination structures in the sense of Definition~\ref{d-appr-eff} as an extension of $\D$-efficient solutions of $f$ in the sense of Definition~\ref{def:sol}(iv), there have been almost {\em no results} for {\em $\ve{k}$-nondominated solutions} as an extension of $\D$-nondominated solutions of $f$ in the sense of Definition~\ref{def:sol}(ii). To the best of our knowledge, the only result in this direction has been obtained by Bao et al. \cite[Theorem~4.7]{best17}, where the conclusions are formulated via a certain auxiliary scalarized function, but not in terms of the given vector-valued mapping. Namely, the scalarized function employed in \cite{best17} is an extended version of the Gerstewitz scalarization function $s : Y\to \mathbb{R}\cup\{\pm \infty\}$ being defined by the formula
\begin{equation}\label{def:s}
s(y) := \operatorname*{inf}\big\{t\in \mathbb{R}\;\big|\; y \in a + tk - \D(y)\big\}.
\end{equation}

To formulate the aforementioned result, we need to recall the following notion, which is also used in the subsequent developments of Section~\ref{new}.

\begin{Definition}[\bf approximate nondominated solutions]\label{d-epsilon-nondominated} Let $f : X \to Y$, $\D : Y \rightrightarrows  Y$, $k\in Y\setminus\{0\}$, and $\ve\geq{0}$. An element $x_{\varepsilon}\in X$ is said to be an {\sc $\ve{k}$-nondominated solution} of $f$ with respect to $\D$ if we have
$$\forall x \in X, \;\; f(x_{\varepsilon})- \ve k \not\in
f(x) + (\D(f(x))\setminus\{0\}).$$
\end{Definition}
For the special case where $\ve =0$, the concept of $\ve{k}$-nondominated solutions of $f$ with respect to $\D$ reduces to $\D$-nondominated solutions of $f$ in the sense of Definition~\ref{def:sol}(ii).\vspace*{0.05in}

The next theorem is an extension of \cite[Theorem~4.7]{best17}, where it is obtained under the assumption that $X$ is a Banach space.

\begin{Theorem}[EVP for approximate nondominated solutions]\label{thm:best} Let $(X,d)$ be a complete metric space while $Y$ is a Banach space, let $k\in Y\setminus\{0\}$, let $\D : Y \rightrightarrows  Y$ be a domination structure, and let $f : X \to Y$ be a vector-valued mapping. Given $\ve\geq{0}$, consider an $\ve{k}$-nondominated solution $x_\ve$ of $f$ with respect to $\D$ and denote $y_\ve := f(x_\ve)$. Impose the following assumptions:
\begin{itemize}
\item[\bf(D1)] {\sc(boundedness condition)} $f$ is bounded from below with respect to the element $\underline{y} \in Y$ and the set $\Theta := \D(\underline{y})$, i.e., $f(x) \in \underline{y} + \Theta$ for all $x\in \dom{f}$. Furthermore,
$\D(\underline{y}) + \D(y) \subseteq \D(y)$ for all $y\in \rge{f}$, $k \in \operatorname*{int}(\Th)$, and $\D(\underline{y})  + \operatorname*{int}(\Th) \subseteq \D(\underline{y})$.

\item[\bf(D2)] {\sc(continuity conditions)} $f$ is continuous over $\dom{f}$, and the domination mapping $\D$ is of closed graph over $\rge{f}$ in the sense that for every sequence of pairs $\{(y_n,v_n)\}$ with $y_n\in f(\O)$ and $v_n\in \D(y_n)$ for all $n\in \mathbb{N}$ the convergence $(y_n,v_n) \to (y_\ast, v_\ast)$ as $n\to\infty$ yields the existence of $x_\ast \in X$ such that $y_\ast = f(x_\ast)$ and $v_\ast \in \D(y_\ast)$.

\item[\bf(D3)] {\sc(scalarization conditions)}
\begin{itemize}
\item[\bf(D3-a)] $\forall y \in \rge{f}$ we have $0 \in \D(y)$ and $\D(y)$ is closed  in $Y$.
\item[\bf(D3-b)] $\forall y\in \rge{f}$ we have $\D(y)+(0,+ \infty)k \subseteq \D(y) \setminus \{0\}$ and $(-\infty,0)k\cap \D(y)=\emptyset$.
\end{itemize}
\end{itemize}
Then there exists an element  $\ox \in\dom{f}$ such that
\begin{itemize}
\item[\bf(i)] $s({f}(\ox)) + \sqrt{{\ve}}d(\ox,x_{\ve}) \leq s({f}(x_\ve))$.

\item[\bf(ii)] $d(\ox,x_{\ve}) \leq \sqrt{\ve}$.

\item[\bf(iii)] $\ox$ is an exact solution of the scalarized function defined by $f_{\ox} := s\circ f + \sqrt{\ve}d(\ox,\cdot)$.
\end{itemize}
\end{Theorem}

\begin{Comment} {\rm It is worth mentioning that the conclusions of Theorem~\ref{thm:best}(i,iii) are formulated in a scalarized form via the scalarization function (\ref{def:s}). It is important to find appropriate assumptions on the given data such that the scalarized function under consideration is bounded from below in condition (D1), and that the continuity condition (D2) is satisfied; cf.\ Lemmas~4.2 and 4.6 in \cite{best17} for more details. Natural questions arise on whether it is possible to weaken the assumptions of this result and/or to obtain conclusions in terms of the nondomination to the given vector-valued mapping by using either a nonlinear scalarization approach, or a nonscalarization approach, or a mixed approach. We develop new results in this direction in the next section.}
\end{Comment}
\vspace*{-0.15in}

%%%%%%%%%%%%%%%%%%%%%%%%%%%%%%%%%%%%%%%%%%%%%%%%%%%%%%%%%%%%%%%%%%%

\section{New Variational Principles in Vector Optimization with Variable Domination Structures in Quasimetric Spaces}\label{new}

This section provides new versions of EVP, which are significantly better than those discussed above and are obtained under weaker assumptions. These new versions of EVP seem to be important for their own sake while having interesting applications to behavioral sciences presented in Section~\ref{behav}. In particular, the results below address $\ve{k}$-efficient solutions of $f$ with respect to $\D$ in the sense of Definition~\ref{d-appr-eff} as well as $\ve{k}$-nondominated solutions of $f$ with respect to $\D$ in the sense of Definition~\ref{d-epsilon-nondominated} under the assumptions that the underlying space $X$ is a quasimetric space (which is essential for applications in Section~\ref{behav}) and that $Y$ is a real linear space. In our approach we use new developments for the Gerstewitz scalarization functions of type \eqref{def:s} given in \cite{gnrt16,w17}.

First we recall the definition of vectorial closedness with respect to a direction and the definition of the Gerstewitz scalarization function. Given a real linear space $Y$ and a nonempty subset $A\subseteq Y$, the {\em vectorial closure} of $A$ in the {\em direction} $k\in Y$ is defined by
$$
\vcl_{k} A :=\big\{y \in Y\;\big|\;\forall \lambda > 0,  \exists t\in [0,\lambda], y + tk \in A\big\}.
$$
We refer the reader to \cite{gnrt16,qh13,w17} for more results and discussions on the directionally vector closedness and its relationships with vector closedness and topological closedness.

\begin{Definition}[\bf nonlinear scalarization functions with domination sets]\label{d-scalar} Let $Y$ be a linear space, let $A$ be a nonempty subset of $Y$, and let $k$ be a nonzero direction in $Y$. The function $\varphi_{A,{k}}: Y \rightarrow \mathbb{R} \cup \{\pm \infty\}$ defined by
\begin{equation}\label{func:ph}
\varphi_{A,k}(y) := \inf\big\{t\in\mathbb{R}\;\big|\;y\in tk-A\big\}\;\mbox{ with }\;\inf{\emptyset} = +\infty
\end{equation}
is called the {\sc Gerstewitz nonlinear scalarization function} generated by the set $A$ and the scalarization direction $k$.
\end{Definition}

By setting $B := A + \mathbb{R}_{+}k$, we get the equalities
$$
\forall y\in Y,\;\varphi_{\vcl_k(B),k}(y) = \varphi_{B,k}(y) = \varphi_{A,k}(y).
$$
\begin{Comment}\label{r-properties}{\rm The the scalarization function $\ph_{A,k}$ was defined in \cite[Theorem~2.3.1]{grtz03} for closed sets $A$ satisfying $A + \mathbb{R}_{+}k \subseteq A$ in real topological vector spaces, where $\ph_{A,k}$ was called a {\em scalarization function with uniform level sets} due to the description of its level sets by
$$
\forall t\in \mathbb{R},\;\operatorname*{Lev}(t;\varphi_{A,k}) = tk - {\vcl}_{k}(B) = tk - {\vcl}_{k}(A + \mathbb{R}_+k).
$$
Furthermore, $\ph_{A,k}$ is {\em translation invariant along the direction} $k$ in the sense that
$$
\forall y\in Y,\;\forall t\in \mathbb{R},\;\ph_{A,k}(y+tk) = t + \ph_{A,k}(y).
$$
Given a subset $B$ of $Y$, $\ph_{A,k}$ is {\em $B$-monotone} in the sense that
$$
a\in b - B \Longrightarrow \ph_{A,k}(a) \leq \ph_{A,k}(b)
$$
if and only if $A+B\subseteq A$. For other properties of the Gerstewitz scalarization functions; see \cite[Theorem~2.3.1]{grtz03}, \cite[Theorem~4]{gnrt16}, and the references therein.}
\end{Comment}

Next we recall some important concepts of quasimetric spaces theory taken from \cite{cob13}.

\begin{Definition}[\bf quasimetric spaces]\label{quasi}
A quasimetric space  is a pair  $(X,q)$ consisting of a set $X$  and  a function $q:X\times X\longmapsto\mathbb R_+:=[0,\infty)$ on $X\times X$ having the following properties:
\begin{itemize}
\item[\bf(i)] $q(x,x^{\prime})\ge 0$ for all $x,x^\prime\in X$ and $q(x,x)=0$ for all $x\in X$ $($positivity$)$.
\item[\bf(ii)] $q(x,x'')\le q(x,x^{\prime})+q(x^{\prime},x'')$ for all $x, x^{\prime},x''\in X$ $($triangle inequality$)$.
\end{itemize}
\end{Definition}

Note that quasimetric spaces may be {\em finite-dimensional}, which is the case of our applications to behavioral science models given in Section~\ref{behav}.

\begin{Definition}[\bf convergence and completeness in quasimetric spaces]\label{quasi1} Let $(X,q)$ be a quasimetric space, and let $\{x_n\}$ be a sequence in $X$.
\begin{itemize}
\item[\bf(i)] The sequence $\{x_n\}$ is said to be {\sc forward-Cauchy} if
for every $\varepsilon >0$ there exists some $N_\varepsilon \in \mathbb{N}$ such that whenever $n\geq N_\varepsilon$ and $m\in \mathbb N$ we have $q(x_n, x_{n+m}) <\varepsilon$.
\item[\bf(ii)]  The sequence $ \{x_n\}$ is said to be {\sc forward-convergent} to $x_\infty$ if  $q(x_n,x_\infty)\to 0$ as $n\to\infty$.

\item[\bf(iii)] The space $(X,q)$ is said to be {\sc forward Hausdorff} if  every forward-convergent sequence has a unique forward-limit.

\item[\bf (iv)] The space $(X,q)$ is {\sc forward-forward-complete} if every forward-Cauchy sequence is forward-convergent.
\end{itemize}
\end{Definition}

Since a quasimetric in {\em not symmetric}, there are the corresponding backward concepts, which can be found in \cite{cob13}.\vspace*{0.05in}

The next concept of generalized Picard sequences of set-valued mappings is taken from \cite{dhm83}.

\begin{Definition}[\bf generalized Picard sequences]\label{pic} A sequence $\{x_n\}$ in a topological space $X$ is called {\sc generalized Picard} for a set-valued mapping $S: X\rightrightarrows  X$  if we have
$$
\forall n\in \N,\;x_{n+1} \in S(x_n).
$$
\end{Definition}

Our approach in this paper is to scalarize a vector-valued mapping $f:X\to{Y}$ by using the Gerstewitz scalarization function $\ph_{A,k}$ defined in (\ref{func:ph}) to construct a generalized Picard sequence of a certain set-valued mapping that converges to the desired element giving us new versions of EVP. Note that the scalarized function $\ph_{A,k}\circ f$ might not be lower semicontinuous.\vspace*{0.05in}

Recall that an Ekeland-type variational principle for set-valued mappings $F$ acting between a complete Hausdorff quasimetric space $X$ and a vector space $Y$ is formulated in \cite[Theorem~4.2]{bms15b}. The binary relation in \cite[Theorem~4.2]{bms15b} is called {\em post-less ordering relation}, which agrees with the efficiency binary relation introduced in Definition~\ref{def:binary}(ii). We derive our first result in this section for {\em $\ve{k}$-efficient solutions} in the sense of Definition~\ref{d-appr-eff} under weaker assumptions than in \cite[Theorem~4.2]{bms15b}, especially those concerning the variable domination structure, monotonicity, and boundedness. It is supposed in \cite[Theorem~4.2]{bms15b} that the set-valued objective mapping $F:X\rightrightarrows Y$ with values in a linear space $Y$ is quasibounded in the following sense.

\begin{Definition}[\bf quasiboundedness]\label{d-quasibounded} A set-valued mapping $F:X\rightrightarrows Y$ with values in a linear space $Y$ is {\sc quasibounded} if there exist a bounded set $M$ in $Y$ such that
$$
\forall x\in X,\;F(x) \subseteq M + \Theta,
$$
where $\Theta$ is the given ordering cone of the image space $Y$.
\end{Definition}

The following major theorem significantly extends the one in \cite[Theorem~4.2]{bms15b} with an essentially different proof. The main tool of our analysis here is a {\em scalarization} technique based on the nonlinear scalarization function introduced in Definition~\ref{d-scalar}.

\begin{Theorem}[\bf variational principle for efficient solutions under variable domination]\label{thm:Pareto}
Let $(X,q)$ be a quasimetric space, let $Y$ be a linear space equipped with a variable domination structure $\D : Y \rightrightarrows Y$, and let $f : X \to Y$ be a vector-valued mapping. Given $k\in Y\setminus\{0\}$, $x_0\in X$, $y_0: = f(x_0)$, $\Th := \D(y_0)$, and $\ve\geq{0}$, we consider the set-valued mapping $W : X\rightrightarrows X$ defined by
\begin{eqnarray}\label{def:W1}
W(x) :=\big\{u\in X\;\big|\;{f(x) - f(u) - \sqrt{\ve}q(x,u)k \in \D(f(x_0))}\big\}
\end{eqnarray}
and the extended-real-valued function $\psi: X \to \R\cup\{\pm\infty\}$ defined by
\begin{equation}\label{f-scal}
\psi(x) := \ph_{\Th,k}(f(x) - f(x_0)) \quad  \big( = \ph_{\Th-f(x_0),k}(f(x))\big),
\end{equation}
where $\ph_{\Th,k}$ is taken from \eqref{func:ph}. Impose the following assumptions:
\begin{itemize}
\item[\bf(E1)] {\sc(boundedness condition)} The function $\psi$ from \eqref{f-scal} is bounded from below over $W(x_0)$.
\item[\bf(E2)] {\sc(limiting monotonicity condition)} For every infinite nonconstant generalized Picard sequence $\{x_n\}$ of the set-valued mapping $W$ from \eqref{def:W1} the convergence of the series\\ $\sum_{n=0}^\infty q(x_n,x_{n+1})$ yields the existence of $x_\ast$ such that
\begin{equation}\label{eq:nonemptyCAP}
\forall n\in \N,\;x_\ast \in W(x_n).
\end{equation}
\item[\bf(E3)] {\sc (scalarization condition)} $\Th$ is $k$-vectorial closed with $0 \in \Th$, $\Theta + \Theta \subseteq \Theta$, $\Theta+  \cone(k) \subseteq \Th$, and $\Th\cap(-\cone(k)) = \{0\}$.
\end{itemize}
Then, there exists $x_\ast \in W(x_0)$ satisfying the inclusion
\begin{equation}\label{W}
W(x_\ast) \subseteq \overline{\{x_\ast\}}:=\big\{u\in X\;\big|\;q(x_\ast,u) = 0\big\}.
\end{equation}
If in addition the condition
\begin{itemize}
\item[\bf(E4)] $(X,q)$ is {\sc forward-Hausdorff}
\end{itemize}
is satisfied, then the conclusions of this theorem reduce to
\begin{itemize}
\item[\bf(i)]  {$f(x_0) - f(x_\ast) - \sqrt{\ve}q(x_0,x_\ast){k} \in \D(f(x_0))$} \; \mbox{ and }

\item[\bf(ii)] $\forall x\in X\setminus\{x_\ast\}$,  $  f(x_\ast) - f(x) -  \sqrt{\ve}q(x_\ast,x)k \not\in \D(f(x_0))$. That is, $x_\ast$ is a $\Theta$-efficient solution of the perturbed function $f_{x_\ast}$, where $f_{x_\ast}: X \to Y$ is defined by
$$
f_{x_\ast}(x) := f(x) + \sqrt{\ve}q(x_\ast,x){k}.
$$
\end{itemize}
Furthermore, imposing the {\sc domination inclusion}
\begin{itemize}
\item[\bf(E5)] $\D(f(x_\ast)) \subseteq \D(f(x_0))$
\end{itemize}
ensures that $x_\ast$ is a $\D$-efficient solution for the perturbed function $f_{x_\ast}$, i.e.,
\begin{equation}\label{D-sol}
\forall x\in X\setminus\{x_\ast\}, f(x_\ast) - f(x) -  \sqrt{\ve}q(x_\ast,x)k \not\in \D(f(x_\ast)).
\end{equation}
If finally the starting point $x_0$ is an {\sc $\ve{k}$-efficient solution} of $f$ with respect to $\D$, then $x_\ast$ can be chosen so that in addition to {\rm(i)} and {\rm(ii)} we have
\begin{itemize}
\item[\bf(iii)] $q(x_0,x_\ast)\leq \sqrt{\ve}$.
\end{itemize}
\end{Theorem}
{\bf Proof}. It is easy to observe from the construction of $W$ in \eqref{def:W1} that
$$
f(x) - f(x_0) \in  -\sqrt{\ve}q(x_0,x)k - \Th \in \R{k} - \Th\;\mbox{ for all }\;u\in W(x_0),
$$
and thus $f(x) - f(x_0) \in \dom\ph_{\Th,k}$. This yields the inclusion $W(x_0) \subseteq \dom{\psi}$ and allows us to construct inductively a generalized Picard sequence satisfying the conclusions of the theorem.

Starting with $x_0$, we assume that $x_n$ is given. Then, choose $x_{n+1} \in W(x_n)$ satisfying
\begin{eqnarray}\label{eqn:choiceTH}
\psi(x_{n+1}) \leq \inf_{u\in W(x_n)} \psi(u) + \frac{1}{2^{n+1}}.
\end{eqnarray}
It is obvious that such an element $x_{n+1}$ exists due to the boundedness from below of the function $\psi$ assumed in (E1). We aim at verifying that the generalized Picard sequence $\{x_n\}$ forward-converges to the desired element by splitting the proof into several steps.\\[1ex]
{\bf Claim~0}: {\em If $u\in W(x)$, then $f(u)\leq_\Th f(x)$ and $W(u) \subseteq W(x)$}. A proof is straightforward. This tells us therefore that
$$
\forall n\in \N,\;{f(x_n) - f(x_{n+1}) \in \Th} \; \mbox{ and } \; W(x_{n+1}) \subseteq W(x_n).
$$
{\bf Claim~1}: {\em For every $u\in W(x_n)$ we have the estimate}
$$
\forall n\in\N, \forall x\in W(x_n),\;\sqrt{\ve}q(x_n,u) \leq \frac{1}{2^{n}}.
$$
Indeed, it follows from $\Theta + \Theta \subseteq \Theta$ assumed in (E3) that the scalarization function $\ph_{\Theta, k}$ is $\Theta$-monotone, i.e., if $v \leq_\Th y$, then $\ph_{\Theta, k}(v) \leq \ph_{\Theta, k}(y)$. Fixing an arbitrary number $n\in \N$ and an arbitrary element $x\in W(x_n)$ yields
\begin{eqnarray*}
x\in W(x_n) & \Longleftrightarrow &   {f(x_n) - f(x) - \sqrt{\ve}q(x_n,x)k \in \Theta}  \\
& \Longleftrightarrow & f(x) - f(x_0) + \sqrt{\ve}q(x_n,x)k \leq_\Th f(x_n) - f(x_0).
\end{eqnarray*}
Since the scalarization function $\ph_{\Th,k}$ given by (\ref{func:ph}) is $\Th$-monotone and translation invariant along the direction $k$ (see Comment~\ref{r-properties}), we have
\begin{eqnarray*}
& & \ph_{\Theta,k}(f(x) - f(x_0) + \sqrt{\ve}q(x_n,x)k) \leq \ph_{\Theta, k}(f(x_{n}) - f(x_0))\\
& \Longrightarrow & \ph_{\Theta,k}(f(x) - f(x_0))  + \sqrt{\ve}q(x_n,x) = \psi(x) + \sqrt{\ve}q(x_n,x) \leq \psi(x_n).
\end{eqnarray*}
This readily implies the inequalities
$$
\sqrt{\ve}q(x_n,x) \leq \psi(x_n) - \psi(x) \leq \psi(x_n) - \inf_{u\in W(x_n)} \psi(u) \leq \psi(x_n) - \inf_{u\in W(x_{n-1})} \psi(u) \leq \frac{1}{2^{n}},
$$
where the last two estimates hold due to $W(x_n)\subseteq W(x_{n-1})$ and (\ref{eqn:choiceTH}), respectively.\\[1ex]
{\bf Claim~2}: {\em The series $\sum_{n=1}^\infty q(x_n,x_{n+1})$ is convergent}. To show this, for every $n\in \N$ we have
\begin{eqnarray}\label{eqn:n-n1}
x_{n+1}\in W(x_n) \Longleftrightarrow  {f(x_n) - f(x_{n+1}) - \sqrt{\ve}q(x_n,x_{n+1})k \in \Th.}
\end{eqnarray}
Summing up these inequalities for $n=0,\ldots, i$ gives us the inclusion
$$
f(x_0) - f(x_{i+1}) -  \sqrt{\ve}(\sum_{n=1}^{i}q(x_n,x_{n+1}))k  \in \Theta,
$$
which can be rewritten in the form
$$
{f(x_{i+1}) - f(x_0) + \sqrt{\ve}(\sum_{n=1}^{i}q(x_n,x_{n+1}))k \in - \Theta.}
$$
Taking into account the $\Th$-monotonicity and transitivity properties of $\ph_{\Th,k}$, we get
\begin{eqnarray*}
& & \ph_{\Theta, k}(f(x_{i+1} - f(x_0) + \sqrt{\ve}(\sum_{n=1}^{i}q(x_n,x_{n+1}))k ) \leq 0\\
& \Longrightarrow & \sqrt{\ve}(\sum_{n=1}^{i}q(x_n,x_{n+1})) \leq - \ph_{\Theta, k}(f(x_{i+1}) - f(x_0) = \psi(x_{i+1}) \leq - \inf_{u\in W(x_0)} \psi(u) <\infty,
\end{eqnarray*}
where the last estimate holds due to (E1). Since $i$ was chosen arbitrary, we arrive at the claimed series convergence
$$
\sqrt{\ve}(\sum_{n=1}^\infty q(x_n,x_{n+1}))<\infty.
$$
{\bf Claim~3}: {\em The inclusion in \eqref{W} is satisfied}. Using the assertion of Claim~2 and the limiting monotonicity condition (E2) ensures the existence of $x_\ast$ satisfying \eqref{eq:nonemptyCAP}. Since we obviously have $x_\ast \in W(x_n) \subseteq W(x_0)$, to get \eqref{W} it is sufficient to prove that
\begin{equation}\label{W1}
\forall u_\ast \in W(x_\ast),\;q(x_\ast,u_\ast) = 0.
\end{equation}
To this end, it is easy to check that
$$
u\in W(x) \Longrightarrow f(x) - f(u) \in \Th  \; \mbox{ and } \; W(u) \subseteq W(x).
$$
Fixing now an arbitrary element $u_\ast \in W(x_\ast)$ gives us by (\ref{eq:nonemptyCAP}) that
$$
\forall n\in \N,\;u_\ast \in W(x_n).
$$
It follows from Claim~1 that $q(x_n,u_\ast) \to 0$ as $n\to\infty$, i.e., $u_\ast$ is a forward-limit of the sequence $\{x_n\}$. Denote further $\overline{\alpha}: = \lim_{n\to\infty}\psi(x_n)$ and show that $\psi(x_\ast) = \overline{\alpha}$. Indeed, the choice of $x_{n+1}$ readily implies that
$$
\psi(x_{n+1}) \leq \inf_{u\in W(x_n)} \psi(u) + \frac{1}{2^{n+1}} \leq \psi(u_\ast) + \frac{1}{2^{n+1}},
$$
where the passage to the limit as $n\to\infty$ yields $\overline{\alpha} \leq \psi(u_\ast)$. On the other hand, we have from the choice of $u_\ast \in W(x_n)$ that
$$
{f(u_\ast) - f(x_0) + \sqrt{\ve}q(x_n,u_{\ast})k \in  f(x_n) - f(x_0) - \Theta.}
$$
Taking into account the $\Th$-monotonicity and the transitivity along the direction $k$ of the scalarization function $\ph_{\Th,k}$ ensures that
$$
\psi(u_\ast) + \sqrt{\ve}q(x_n,u_{\ast}) \leq   \psi(x_n),
$$
and thus we get by passing to the limit as $n\to\infty$ that $\psi(u_\ast) \leq \overline{\alpha}$. Hence $\psi(u_\ast) = \overline{\alpha}$. Since $x_\ast \in W(x_\ast)$, we have $\psi(x_\ast) = \overline{\alpha}$ as claimed. It now follows from $u_\ast\in W(x_\ast)$ that
$$
{f(u_\ast) -f(x_0) + \sqrt{\ve}q(x_\ast,u_{\ast})k \in  f(x_\ast) - f(x_0) - \Theta.}
$$
The aforementioned properties of the scalarization function $\ph_{\Th,k}$ lead us to
$$
\psi(u_\ast) + \sqrt{\ve}q(x_\ast,u_{\ast}) \leq   \psi(x_\ast).
$$
Substituting  $\psi(u_\ast) = \psi(x_\ast) = \overline{\alpha}$ into the last inequality, we have $q(x_\ast,u_{\ast}) \leq 0$ and hence $q(x_\ast,u_{\ast}) = 0$, which verifies \eqref{W1}
Since $u_\ast$ was chosen arbitrarily in $W(x_\ast)$, we have
\begin{eqnarray}
\label{eqn:conclusion}
x_\ast \in W(x_0) \; \mbox{ and } \; W(x_\ast) \subseteq \overline{\{x_\ast\}}.
\end{eqnarray}
This ensures the fulfillment of \eqref{W} and thus completes the proof of Claim~3.\\[1ex]
{\bf Claim~4}: {\em Imposing {\rm(E4)} gives us assertions {\rm(i)} and {\rm(ii)} of the theorem}. Assumption (E4) tells us that the forward-limit is unique if exists. Thus $\overline{\{x_\ast\}} = \{x_\ast\}$, and the two inclusions in (\ref{eqn:conclusion}) reduce to assertions (i) and (ii), respectively, by the construction of the sets $W(x)$ in \eqref{def:W1}.\\[1ex]
{\bf Claim~5}: {\em The domination inclusion in {\rm(E5)} yields \eqref{D-sol}.} Arguing by contraposition, suppose that $x_\ast$ is not a $\D$-efficient solution of the perturbed function $f_{x_\ast}$. Then, we find $x \neq x_\ast$ such that
$$
f(x) + \sqrt{\ve}q(x_\ast,x)k \in f(x_\ast) - \D(x_\ast) \st{{ (E5)}}{\subseteq} f(x_\ast) - \Th,
$$
which clearly contradicts (ii).\\[1ex]
{\bf Claim~6}: {\em If $x_0$ is an $\ve{k}$-efficient solution of $f$ with respect to $\D$, then we have {\rm(iii)}.} We again argue by contraposition and suppose that (iii) fails, i.e., $q(x_0,x_\ast) > \sqrt{\ve}$. Then, (i) yields
\begin{eqnarray*}
f(x_\ast) & \in &  f(x_0) - \sqrt{\ve}q(x_0,x_\ast)k - \D(f(x_0)) \\
& = &  f(x_0) - \ve{k} -\sqrt{\ve}(q(x_0,x_\ve) - \sqrt{\ve})k - \D(f(x_0)) \\
& \st{(E2)}{\subseteq} &  f(x_0) - \ve{k} - \D(f(x_0)).
\end{eqnarray*}
This readily contradicts the assumption that $x_0$ is an $\ve{k}$-efficient solution of $f$ with respect to $\D$ and thus completes the proof of the theorem. $\h$\vspace*{0.05in}

It is clear that the obtained Theorem~\ref{thm:Pareto} weakens and/or drops many assumptions in Theorem~\ref{thm:bms15} and Theorem~\ref{thm:benahn}. Let us comment on the major assumptions imposed in Theorem~\ref{thm:Pareto}.

\begin{Comment} {\bf (a)} {\rm Condition (E1) allows us to extend EVP to the class of vector-valued mappings having their image spaces as arbitrary linear spaces. When the image space happens to be a normed one, (E1) is weaker than the quasiboundedness condition (A1) and the boundedness condition (B1). In fact, (E1) can be equivalently written as
$$
\exists e\in Y,\;\forall x\in W(x_0), f(x) - f(x_0) \in e + Y\setminus(-\Theta).
$$
{\bf(b)} Condition (E2) is better than conditions (A2) and (B2) as shown in the next two propositions.}
\end{Comment}

The following simple example shows that the boundedness from below condition (E1) imposed in Theorem~\ref{thm:Pareto} is weaker than the quasiboundedness assumption (A1) imposed in Theorem~\ref{thm:bms15}.

\begin{Example}{\bf(boundedness from below condition of Theorem~\ref{thm:Pareto} versus quasiboundedness)}\label{exa-bound} {\rm Consider a vector-valued mapping $f: \R \rightrightarrows  \R^2$ defined by
\begin{eqnarray*}
f(x) := \left\{\begin{array}{lll}
(0,-x) &\mbox{if}& x\geq 0,\\
(x,0) &\mbox{if}&x<0,
\end{array}\right.
\end{eqnarray*}
and consider a fixed domination structure $\D: \R^2 \rightrightarrows  \R^2$ with $\D(y) \equiv \R^2_+$. In this case we have $\Th = \R^2_+$. Take $k:=(1,1)\in \Th\setminus -\Th$ and $q(x,u):= |x-u|$. Then, $W(x) = \{x\}$ for all $x\in \R$. The mapping $f$ clearly satisfies the boundedness from below condition (E1), but it is not quasibounded in the sense of (A1) as required in Theorem~\ref{thm:bms15}. Indeed, we have
$$
\rge{f} = \cone\{(0,-1),(-1,0)\}.
$$}
\end{Example}

Next we present an illustrative example for other assumptions of Theorem~\ref{thm:Pareto}.

\begin{Example}[\bf illustrating the assumptions of Theorem~\ref{thm:Pareto}]\label{exa2} {\rm Let $X := \R$ and $Y := \R^2$, and let $f: X\to Y$ be defined by
$$
f(x) := \begin{cases}
(x, 2^{x} - 1) & \mbox{ if }\; x< 0,\\
 (x,1)  & \mbox{ if }\; x\geq 0.
\end{cases}
$$
The domination structure $\D: \R^2 \rightrightarrows  \R^2$ is given by
$$
\D(y) := \begin{cases}
\conv\cone\big\{(1,0), (|y_1|,|y_2|)\big\} & \mbox{ if }\; y_1< 0 \; \mbox{ and } \; y_2 < 0,\\
\R^2_+  & \mbox{ otherwise. }
\end{cases}
$$
Take $\ve = 1$, $x_0 = 0$, $k = (1,1)$, and $d(x,u) = 1/2|x-u|$. In this case we have $\psi(x) = \ph_{\R^2_+,k}$. It is easy to check that:
\begin{itemize}
\item[\bf a)] $W(0) = [-0.5,0]$ and $W(-0.5) = \{-0.5\}$.
\item[\bf b)] $f$ is $\ph_{\R^2_+,k}$-bounded from below.
\item[\bf c)] $f$ is not $\R^2_+$-lower semicontinuous since
$$
\operatorname*{lev}(f,0_{\R^2}) =\big\{x\in X\;\big|\;f(x) \in 0 - \R^2_+\big\} = (-\infty,0)
$$
is not a closed set in $\R$.
\item[\bf d)] $f(-0.5) =(-0.5,-1 + 1/\sqrt{2})$, $f(0) = (0,1)$, $\D(f(-0.5)) = \conv\cone\{(1,0), (0.5, 1- 1/\sqrt{2})\}$, and $\D(f(0)) = \R^2_+$. It is obvious that $f(-0.5) \leq_{\D(f(x_0))} f(0)$ and $\D(f(-0.5)) \subseteq \D(f(0))$, and hence condition (E5) is satisfied.
\item[\bf e)] Condition (E2) holds since for any nonconstant generalized Picard sequences in $W_0$ (without loss of generality it can assumed that $x_n < 0$) we have that $W(x_n)$ is closed whenever $n\in\N$. Then, the existence of $x_\ast$ follows from the classical Cantor theorem.
\end{itemize}}
\end{Example}

Our next goal is to derive efficient conditions expressed entirely via the given problem data ensuring the fulfillment of the limiting monotonicity assumption (E2) of Theorem~\ref{thm:Pareto}.
First we introduce the following new notion.

\begin{Definition}{\bf(decreasing lower semicontinuity with respect to domination sets)}\label{dec-lsc} A vector-valued mapping $f: X\to Y$ is said to be {\sc $\Th$-decreasing lower semicontinuous} over a set $A$ if for every forward-convergent sequence $\{x_n\} \subseteq A$ with a forward-limit $x_\ast$, the $\Th$-decreasing monotonicity of $\{f(x_n)\}$ with respect to the {\sc domination set} $\Th$ (i.e., $f(x_{n+1}) \leq_\Th f(x_n)$ for all $n\in \N$) implies that
$$
\forall n\in \N,\;f(x_\ast) \leq_\Th f(x_n).
$$
\end{Definition}

\begin{Proposition} [\bf sufficient conditions for limiting monotonicity]\label{sufE2a} Let $(X,q)$ be a forward-complete quasimetric space, and let the vector-valued mapping $f: X\to Y$ be $\Th$-decreasing lower semicontinuous over $W(x_0)$. Then the set-valued mapping $W$ from \eqref{W} satisfies the limiting monotonicity condition {\rm(E2)} of Theorem~{\rm\ref{thm:Pareto}}.
\end{Proposition}
{\bf Proof}. Take an arbitrary generalized Picard sequence $\{x_n\}\subseteq W(x_0)$ of the set-valued mapping $W$ satisfying $\sum_{n=0}^\infty q(x_n,x_{n+1}) = \ell <\infty$. Then, for each $\ve > 0$ there exists $N_\ve\in \N$ such that
$$
\sum_{n=0}^{N_\ve-1} q(x_n,x_{n+1}) \geq \ell - \ve \; \mbox{ and } \;\sum_{n=N_\ve}^\infty q(x_n,x_{n+1}) \leq \ve.
$$
For every $i, j \geq N_\ve$ with $j > i$ we have
$$
q(x_i,x_j) \leq \sum_{n=i}^{j-1} q(x_n,x_{n+1}) \leq \sum_{n=N_\ve}^\infty q(x_n,x_{n+1}) \leq \ve.
$$
Hence $\{x_n\}$ is a forward-Cauchy sequence. Since $(X,q)$ is forward-forward-complete, it forward-converges to some forward-limit $x_\ast \in X$. Thus $\{x_n\}$ is a generalized Picard sequence of $W$, i.e.,
$$
\forall n\in \N,\;x_{n+1} \in W(x_n),
$$
which implies that $f(x_{n+1}) \leq_\Th f(x_n)$ for all $n\in \N$ by (E4). The imposed $\Th$-decreasing lower semicontinuity of $f$ ensures that
$$
\forall n\in \N,\,f(x_{\ast}) \leq_\Th f(x_n).
$$
Since $\{x_n\}$ is a generalized Picard sequence of $W$, it follows that
\begin{eqnarray*}
\forall n\in \N,\;f(x_{n+1}) + \sqrt{\ve}q(x_n,x_{n+1})k \leq_\Th  f(x_n).
\end{eqnarray*}
Summing up these relations from $n=i$ to $i+j$ while taking into account (E4) and the triangle inequality for the quasimetric, we have
$$
f(x_{i+j}) + \sqrt{\ve}q(x_i,x_{i+j})k \leq_\Th f(x_i).
$$
Adding this to $f(x_\ast) \leq_\Th f(x_{i+j})$ yields
$$
f(x_{\ast}) + \sqrt{\ve}q(x_i,x_{i+j})k \leq_\Th f(x_i).
$$
Taking again into account the triangle inequality for the quasimetric and (E4) gives us the estimates
\begin{eqnarray*}
f(x_\ast) + \sqrt{\ve}q(x_i,x_\ast)k + \sqrt{q}(x_{i+j}, x_\ast)k
\leq_\Th f(x_{\ast}) + \sqrt{\ve}q(x_i,x_{i+j})k \leq_\Th f(x_i).
\end{eqnarray*}
Since $j$ was chosen arbitrarily, $\Th$ is $k$-vectorial closed, and since $\lim_{n\to\infty} q(x_n,x_\ast) = 0$, the passage above to the limit as $j\to\infty$ yields
$$
f(x_{\ast}) - \sqrt{\ve}q(x_i,x_\ast)k  \in f(x_i) - \Th,
$$
i.e., $x_\ast \in W(x_i)$. The latter verifies the fulfillment of (E2), since $i$ was also chosen arbitrarily. $\h$

\begin{Proposition}[\bf other sufficient conditions for the fulfillment of (E2)]\label{sufE2b} Let $\Th$ satisfy condition {\rm(E3)}, and let $f$ be $(\Th,k)$-lower semicontinuous in the sense of  Soleimani {\rm\cite{s14}}:
$$
\forall x\in W(x_0),\forall t\in \R,\;L(x,t) := \big\{u\in X\;\big|\; f(u) \in f(x) + tk - \Theta\big\} \; \mbox{ is closed in } \; X.
$$
Then, $f$ satisfies condition {\rm(E2)} of Theorem~{\rm\ref{thm:Pareto}}.
\end{Proposition}
{\bf Proof}. Pick an arbitrary generalized Picard sequence $\{x_n\}\subseteq W(x_0)$ of the set-valued mapping $W$ satisfying $\sum_{n=0}^\infty q(x_n,x_{n+1}) = \ell <\infty$, and then fix an arbitrary number $n\in \N$. For every $i\in \N$, it is not difficult to check that
$$
f(x_{n+i}) + \sqrt{\ve}q(x_{n},x_{n+i})k \leq_\Th  f(x_n)
$$
by using (E3) and the triangle inequality for the quasimetric. We can further proceed as follow:
\begin{eqnarray*}
f(x_{n+i}) & \in &  f(x_n) - \sqrt{\ve}q(x_{n},x_{n+i})k - \Th \\
& = & f(x_n) - \sqrt{\ve}q(x_n,x_\ast) + \sqrt{\ve}q(x_{n+i},x_\ast)k \\
& & \hspace*{1.5in}- \sqrt{\ve}\big(q(x_{n+i},x_\ast)k - q(x_n,x_\ast)+ q(x_n,x_{n+i})\big)k - \Th\\
& \subseteq & f(x_n) - \sqrt{\ve}q(x_n,x_\ast)k + \sqrt{\ve}q(x_{n+i},x_\ast)k -\Th.
\end{eqnarray*}
Since $\lim_{i\to\infty} q(x_{n+i},x_\ast) = 0$, for every $\delta>0$ there is $N_\delta \in \N$ such that
$$
f(x_{n+i}) \in f(x_n) - \sqrt{\ve}q(x_n,x_\ast)k + \sqrt{\ve}(\delta)k - \Th.
$$
Since $f$ is $(\Theta,k)$-lower semicontinuous, we have
$$
f(x_{\ast}) \in f(x_n) - \sqrt{\ve}q(x_n,x_\ast)k + \sqrt{\ve}\delta{k} - \Th.
$$
Remembering that $\delta>0$ was chosen arbitrarily while $\Th$ is $k$-vectorially closed, it follows that
$$
f(x_\ast) \in f(x_0) + f(x_{n}) - \sqrt{\ve}q(x_n,x_\ast)k - \Th,
$$
i.e., $x_\ast \in W(x_n)$. Since $n$ was  also chosen arbitrarily,  condition (E2) is verified. $\h$\vspace*{0.15in}

Our second major result is a new variational principle for $\ve{k}$-{\em nondominated solutions} of $f$ with respect to $\D$ in the sense of Definition~\ref{d-epsilon-nondominated}. Note to this end that in \cite[Theorem~4.5]{bms15b} an Ekeland-type variational principle is obtained for set-valued objective mappings $F:X\rightrightarrows  Y$ between a complete Hausdorff quasimetric space $X$ and a linear space $Y$. The binary relation considered in \cite[Theorem~4.5]{bms15b} is called {\em pre-less ordering relation}, which agrees with the nondomination binary relation introduced in Definition~\ref{def:binary}(i).\vspace*{0.05in}

Now we establish our variational principle in quasimetric spaces for {\em nondominated solutions} under weaker assumptions than in \cite[Theorem~4.5]{bms15b}, especially those concerning variable domination structures, monotonicity and boundedness. Furthermore, as in the case of Theorem~\ref{thm:Pareto}, the proof of the following theorem is based on the nonlinear scalarization function from Definition \ref{d-scalar}, which is  significantly different from the proof of \cite[Theorem~4.5]{bms15b}.

\begin{Theorem}{\bf(variational principle for nondominated solutions under variable domination)}\label{thm:EVP-dominated} Let $(X,q)$ be a quasimetric space, let $Y$ be a linear space, let $\D:Y\rightrightarrows  Y$ be a domination structure on $Y$ with the nondomination relation $\leq_{N}$, let $k\in Y\setminus\{0\}$, and let $\Theta:= \D(f(x_0))$. Given $x_0 \in X$ and $\ve\geq{0}$, define the set-valued mapping $W = W_{f,q,\ve}: X\rightrightarrows X$ by
\begin{eqnarray}\label{def:W}
W(x) :=\big\{u\in X \;\big|\;{f(x) -\sqrt{\ve}q(x,u)k - f(u) \in \D(f(u))}\big\};
\end{eqnarray}
where we drop the parameters $f, q$, and $\ve$ in the notation of $W$ for simplicity. Consider also the extended-real-valued function $\psi: X \to \R\cup\{\pm\infty\}$ given by
\begin{eqnarray}\label{def:s1}
\psi(x): = \ph_{\Th,k}(f(x) - f(x_0)) \quad  \big( = \ph_{\Th-f(x_0),k}(f(x))\big),
\end{eqnarray}
where $\ph_{\Th,k}$ was introduced in \eqref{func:ph}. Impose the following assumptions:
\begin{itemize}
\item[\bf(F1)] {\sc(boundedness condition)} $\psi$ is bounded from below on $W(x_0)$; i.e., there exists $\tau\in \R$ such that $\psi(x)\geq \tau$ for all $x\in W(x_0)$.

\item[\bf(F2)] {\sc(limiting monotonicity condition)} for every generalized Picard sequence $\{x_n\}$ of the set-valued mapping $W$, the convergence of the series $\sum_{n=0}^\infty q(x_n,x_{n+1})$ yields the existence of $x_\ast$ satisfying $x_\ast \in W(x_n)$ for all $n\in \N$.

\item[\bf(F3)] {\sc(scalarization conditions)}
\begin{itemize}
\item[\bf(F3-a)] $\Th$ is $k$-vectorially closed, $\Theta + \Theta \subseteq \Theta$, $\Theta+  \cone(k) \subseteq \Th$, and $\Th \cap -\cone(k) = \{0\}$.
\item[\bf(F3-b)] $\forall x\in W(x_0)$, $\D(f(x)) + \cone(k) \subseteq \D(f(x))$.
\item[\bf(F3-c)] $\forall {f(u) - f(w) \in \D(f(w))}$, $\D(f(w)) + \D(f(u)) \subseteq \D(f(w))$.
\item[\bf(F3-d)] $\Theta^u\subseteq \Th$, where $\Theta^u := \cup\{\D(f(x)) \;:\; x\in W(x_0)\}$.
\end{itemize}
\end{itemize}
Then, there exists $x_\ast\in W(x_0)$ such that $W(x_\ast) \subseteq \overline{\{x_\ast\}}$, where
$$
\overline{\{x_\ast\}} :=\big\{u\in X \;\big|\; q(x_\ast,u) = 0\big\}.
$$
Assuming furthermore that
\begin{itemize}
\item[\bf(F4)] the forward-limit of a forward-convergent sequence in the quasimetric space $(X,q)$ is {\sc unique},
\end{itemize}
the conclusions of the theorem can be written in the equivalent form
\begin{itemize}
\item[\bf(i)] $f(x_0) - \sqrt{\ve}q(x_0,x_\ast) - f(x_\ast)  \in \D(f(x_\ast))$,
\item[\bf(ii)] $\forall x \neq x_\ast,\; f(x_\ast) - \sqrt{\ve}q(x_\ast,x)k -f(x) \not\in \D(f(x))$.
\end{itemize}
If finally the starting point $x_{0}$ is an {\sc $\ve{k}$--nondominated solution} of $f$ with respect to $\D$, i.e.,
\begin{eqnarray}\label{appr}
\forall x\in X,\;f(x_0) - \ve{k} - f(x)\not\in \D(f(x)),
\end{eqnarray}
then we have in addition to {\rm(i)} and {\rm(ii)} that $x_\ast$ satisfies the localization condition
\begin{itemize}
\item[\bf(iii)] $q(x_{0},x_\ast)\leq \sqrt{\ve}$.
\end{itemize}
\end{Theorem}
{\bf Proof}. Let us construct inductively a generalized Picard sequence that satisfies all the requirements of the theorem. Observe that for every $u\in W(x_0)$ we have
$$
f(u) \in f(x_0) - \sqrt{\ve}q(x_0,u)k - \D(f(u)) \st{(F3-a)}{\subseteq} f(x_0) - \D(f(u)) \st{(F3-d)}{\subseteq} f(x_0) - \Theta,
$$
which clearly implies that $\ph_{\Theta, k}(f(u)-f(x_0))$ is finite due to (F1). This gives us $W(x_0) \subseteq \dom{\psi}$.

Starting with $x_0$, suppose that $x_n$ is defined, and then choose $x_{n+1} \in W(x_n)$ such that
\begin{eqnarray}
 \psi(x_{n+1}) \leq \inf_{u\in W(x_n)}  \psi(u) + \frac{1}{2^{n+1}}.
\end{eqnarray}
{\bf Claim~1}: {\em If $u\in W(x)$, then $W(u) \subseteq W(x)$}. To verify this, fix an arbitrary element $w\in W(u)$ and by using the construction of $W$ get the inclusions
\begin{eqnarray}
\label{eqn:comp}
 {f(x) - \sqrt{\ve}q(x,u)k - f(u) \in  \D(f(u)) \; \mbox{ and } f(u) -\sqrt{\ve}q(u,w)k - f(w) \in \D(f(w)).}
\end{eqnarray}
By (F3-b), the second inclusion in (\ref{eqn:comp}) yields  {$ f(u) - f(w) \in \D(f(w))$, and thus $\D(f(u)) + \D(f(w)) \subseteq \D(f(w))$ due to (F3-c).} Combining both inclusions in (\ref{eqn:comp}) tells us that
$$
{  f(x) -\sqrt{\ve}q(x,u)k - \sqrt{\ve}q(u,w)k - f(w)  \in \D(f(w)).}
$$
Employing (F3-b) and the triangle inequality of the quasimetric $q$, we obtain
\begin{eqnarray*}
& & f(x) -\sqrt{\ve}q(x,u)k - \sqrt{\ve}q(u,w)k - f(w)   \in   \D(f(w))\\
& \Longleftrightarrow & f(x) -\sqrt{\ve}q(x,w)k - f(w)   \in   \D(f(w)) + \sqrt{\ve}(q(x,u) + q(u,w) - q(x,w))k\\
& \Longrightarrow & f(x) -\sqrt{\ve}q(x,w)k - f(w)   \in   \D(f(w)),
\end{eqnarray*}
which readily implies that $w\in W(x)$. Since $w$ was chosen arbitrarily in $W(u)$, it gives us
$W(u) \subseteq W(x)$ as asserted in Claim~1.\\[1ex]
{\bf Claim~2}: {\em For all $n\in \N$ and for every $u\in W(x_n)$ we have $\sqrt{\ve}q(x_n,u)\leq \frac{1}{2^{n}}$}. Picking an arbitrary element $u\in W(x_n)$, observe that
{\begin{eqnarray*}
& &   f(x_n) -\sqrt{\ve}q(x_n,u)k - f(u) \in \D(f(u))\\
&\st{(F3-d)}{\Longrightarrow} &  f(x_n) -\sqrt{\ve}q(x_n,u)k -f(u) \in \Theta\\
& \Longleftrightarrow & f(u) - f(x_0) \in f(x_n) - f(x_0) -\sqrt{\ve}q(x_1,x_{n+1})k - \Theta.
\end{eqnarray*}}
Since $\Theta +\Theta \subseteq \Theta$ by (F3-a), the scalarization function $\ph_{\Theta, k}$ is $\Theta$-monotone. This leads us to
\begin{eqnarray*}
\ph_{\Theta, k}(f(u) - f(x_0)) & \leq & \ph_{\Theta, k}(f(x_{n}) - f(x_0) - \sqrt{\ve}q(x_n,u)k)\\
& = & \ph_{\Theta, k}(f(x_{n}) - f(x_0)) - \sqrt{\ve}q(x_n,u).
\end{eqnarray*}
Using the assertion of Claim~1, we have $W(x_{n})\subseteq W(x_{n-1})$ and thus arrive at the estimates
$$
\sqrt{\ve}q(x_n,u) \leq \psi(x_n) -  \psi(u) \leq  \psi(x_n) - \inf_{u\in W(x_n)} \psi(u) \leq  \psi(x_n) - \inf_{u\in W(x_{n-1})}  \psi(u) \leq \frac{1}{2^{n}},
$$
where the two last inequalities hold due to $W(x_n)\subseteq W(x_{n-1})$ and (\ref{eqn:choiceTH}). This verifies the claim.\\[1ex]
{\bf Claim~3}: {\em The series $\sum_{n=1}^\infty q(x_n,x_{n+1})$ is convergent}. For every $n\in \N$ we have $x_{n+1}\in W(x_n)$, which is equivalent to the inclusion
\begin{eqnarray}\label{eqn:choiceW}
{ f(x_n) -\sqrt{\ve}q(x_n,x_{n+1})k - f(x_{n+1}) \in \D(f(x_{n+1})).}
\end{eqnarray}
It follows from (F2-b) that {$f(x_n) - f(x_{n+1}) \in\D(f(x_{n+1}))$}. Then, using (F3-c) gives us $\D(f(x_{n+1})) + \D(f(x_n)) \subseteq \D(f(x_{n+1}))$. Summing up the inclusions in (\ref{eqn:choiceW}) for $n=0,\ldots, i$, we get that
$$
{f(x_0) - \sqrt{\ve}(\sum_{n=1}^{i}q(x_n,x_{n+1}))k - f(x_{i+1}) \in \D(f(x_{i+1})).}
$$
Taking further (F3-d) into account ensures that
$$
{f(x_{i+1}) - f(x_0) + \sqrt{\ve}(\sum_{n=1}^{i}q(x_n,x_{n+1}))k \in - \Theta.}
$$
Since $\ph_{\Th,k}$ is $\Th$-monotone and translation invariant along the direction $k$ as discussed in Comment~\ref{r-properties}, we obtain the inequality
$$
\ph_{\Theta, k}(f(x_{i+1}) - f(x_0)) \leq -\sqrt{\ve}(\sum_{n=1}^{i})q(x_n,x_{n+1})) .
$$
Then, the boundedness condition (F1) tells us that
$$
\sqrt{\ve}(\sum_{n=1}^{i}q(x_n,x_{n+1}))  \leq - \psi(x_{i+1}) \leq  -\inf_{u\in W(x_0)}  \psi(u) <\infty.
$$
Since $i$ was chosen arbitrarily, we arrive at
$$
\sqrt{\ve}(\sum_{n=1}^\infty q(x_n,x_{n+1}))<\infty,
$$
which readily verifies this claim.\\[1ex]
{\bf Claim~4}: {\em We have the inclusion $W(x_\ast) \subseteq \overline{\{x_\ast\}}$}. Since (F2) yields the existence of $x_\ast$ with $x_\ast \in W(x_n)$ for all $n\in \N$, to justify this claim it is sufficient to show that
$$
\forall u_\ast\in W(x_\ast),\;q(x_\ast,u_\ast) = 0.
$$
Fix an arbitrary element $u_\ast\in W(x_\ast)$. It follows from Claim~3 with $u_\ast \in W(x_n)$ and the assertion of Claim~2 that $\sqrt{\ve}q(x_n,u_\ast)\leq \frac{1}{2^{n}}$ for all $n\in \N$, i.e., $u_\ast$ is a forward-limit of the sequence $\{x_n\}$.

Denoting $\overline{\alpha}:= \lim_{n\to\infty} \psi(x_n)$, we intend to prove that $\psi(u_\ast) = \overline{\alpha}$. Indeed, it follows from the choice of $x_{n+1}$ that
$$
\psi(x_{n+1}) \leq \inf_{u\in W(x_n)} \psi(u) - \frac{1}{2^{n+1}} \leq \psi(u_\ast) -  \frac{1}{2^{n+1}},
$$
where the passage to the limit as $n\to\infty$ yields $\overline{\alpha} \leq \psi(u_\ast)$. To verify the opposite inequality, deduce from $u_\ast \in W(x_n)$ that
\begin{eqnarray*}
f(u_\ast) - f(x_0) & \in &   f(x_n) - f(x_0) - \sqrt{\ve}q(x_n,u_{\ast})k - \D(f(u_\ast))\\
& \st{(F3-d)}{\subseteq} & f(x_n) - f(x_0) - \sqrt{\ve}q(x_n,u_{\ast})k - \Th.
\end{eqnarray*}
Since the scalarization function $\ph_{\Th,k}$ is $\Th$-monotone and translation invariant along the direction $k$ (see Comment~\ref{r-properties}), we get
$$
\psi(u_\ast) \leq   \psi(x_n) - \sqrt{\ve}q(x_n,u_{\ast}).
$$
Passing there to the limit as $n\to\infty$ yields $\psi(u_\ast) \leq \overline{\alpha}$, and thus $\psi(u_\ast) = \overline{\alpha}$. Furthermore, we obtain $\psi(x_\ast) = \overline{\alpha}$ since $x_\ast \in W(x_\ast)$.

Deduce now from $u_\ast\in W(x_\ast)$ the inclusion
$$
f(u_\ast) -f(x_0) \in   f(x_\ast) - f(x_0) - \sqrt{\ve}q(x_\ast,u_{\ast})k - \Theta.
$$
Employing again the aforementioned properties of the scalarization function $\ph_{\Th,k}$ implies that
$$
\psi(u_\ast) + \sqrt{\ve}q(x_\ast,u_{\ast}) \leq    \psi(x_\ast).
$$
Substituting $ \psi(u_\ast) =  \psi(x_\ast) = \overline{\alpha}$ into the last inequality, we arrive at $q(x_\ast,u_{\ast}) \leq 0$ and thus justify that $q(x_\ast,u_{\ast}) = 0$. This readily verifies Claim~4 due to the construction of the set $W(x_\ast)$ in \eqref{def:W}.\\[1ex]
{\bf Claim~5}: {\em Assertions~{\rm(i)} and {\rm(ii)} under assumption {\rm(F4)}}. Assuming (F4) ensures that $\overline{\{x_\ast\}} = \{x_\ast\}$. This yields by the constructions above that $x_\ast\in W(x_0)$ and $W(x_\ast) = \{x_\ast\}$, which can be equivalently written as forms (i) and (ii) of the theorem.\\[1ex]
{\bf Claim~6}: {\em Completing the proof}. It remains to estimate the quasidistance between $x_0$ and $x_\ve$ when $x_0$ is an $\ve{k}$--nondominated solution of $f$ with respect to $\D$. Indeed, if (iii) fails, we have $q(x_0,x_\ve) > \sqrt{\ve}$. Then, it follows from (i) and (F3-b) that
\begin{eqnarray*}
f(x_\ve) & \in &  f(x_0) - \sqrt{\ve}q(x_0,x_\ve) - \D(f(x_\ve)) = f(x_0) - \ve{k} -\sqrt{\ve}(q(x_0,x_\ve) - \sqrt{\ve})k - \D(f(x_\ve)) \\
& \subseteq &  f(x_0) - \ve{k} - \D(f(x_\ve)),
\end{eqnarray*}
which clearly contradicts the choice of $x_0$ and thus completes the proof of the theorem. $\h$

\begin{Comment} {\bf(a)} {\rm In the proof of Theorem~\ref{thm:EVP-dominated} we provide a new way to construct a generalized Picard sequence of $W$ by using the scalarization function $\ph_{\Th,k}$. This seems to be more effective than the procedure in \cite[Theorem~3.4]{bm10}: find $x_{n+1} \in W(x_n)$ such that
$$
q(x_n,x_{n+1}) \geq \sup_{u\in W(x_n)} q(x_n,u)-\frac{1}{2^{n+1}}.
$$

{\bf(b)} The boundedness condition (F1) is less restrictive than the requirement that $f$ is quasibounded with respect to $\Theta$ (see Definition~\ref{d-quasibounded}) as supposed in (H3) of \cite[Theorem~4.5]{bms15b}.

{\bf(c)} The limiting monotonicity condition (F2) is related to the level-decreasing-closedness condition (H4') used in \cite[Theorem~4.5]{bms15b} under the name of ``pre-less preorder."

{\bf(d)} Since there is only one domination set $\Th = \D(f(x_0))$ satisfying condition (F3-a), the additional condition (F3-d) is essential in comparison with (D3). Condition (F3-c) guarantees that the domination binary relation $\leq_N$ is transitive.}
\end{Comment}

The next proposition provides efficient conditions in terms of the given data of the problem that ensure the fulfillment of the boundedness assumption (F1) in Theorem~\ref{thm:EVP-dominated} expressed therein via the auxiliary function \eqref{def:s1}.

\begin{Proposition}{\bf(sufficient conditions for the boundedness assumption of Theorem~\ref{thm:EVP-dominated})}\label{suffF1} Let $Y$ be a normed space, and let a set $\Th$ be topologically closed in $Y$  with $k\in \Th\setminus (-\Theta)$. Assume in addition that {\rm(F3-a)} holds, and that the mapping $f$ in Theorem~{\rm\ref{thm:EVP-dominated}} is quasibounded over $W(x_0)$. Then, the boundedness condition {\rm(F1)} is satisfied.
\end{Proposition}
{\bf Proof}. Since $f$ is quasibounded over $W(x_0)$, there exists a bounded set $M$ such that
\begin{eqnarray}
\label{eqn:boundedness}
\forall x\in W(x_0),\;f(x) - f(x_0) \in M + \Theta.
\end{eqnarray}
Arguing by contraposition, assume that the function $\psi$ from (\ref{def:s1}) is not bounded from below over $W(x_0)$. Then, there is a sequence $\{x_n\} \subseteq W(x_0)$ such that $\psi(x_n) \leq -n$, and thus
$$
\forall n\in \N,\;f(x_n) - f(x_0) \in -nk - \Theta.
$$
Fixing an arbitrary number $n\in \N$, we find $\theta_k \in \Theta$ such that $f(x_n) - f(x_0) = -nk - \theta_k$. It follows from (\ref{eqn:boundedness}) that $-nk - \theta_k \in M + \Theta$, and then by (F3-a) we have
$$
k \in \frac{1}{n}M - \Theta.
$$
Since $n$ was chosen arbitrarily in $\N$ while $\Th$ is a closed set, the passage to limit as $n\to\infty$ yields $k\in -\Theta$. This clearly contradicts the choice of $k\in \Th\setminus (-\Theta)$ and thus completes the proof. $\h$\vspace*{0.05in}

Now we present an example illustrating the fulfillment of all the scalarization conditions in assumption (F3) of Theorem~\ref{thm:EVP-dominated}.

\begin{Example}[\bf scalarization conditions of Theorem~\ref{thm:EVP-dominated}]\label{exa3} {\rm In the setting of $Y = \R^2$ with $k = (1,1)$, consider the ordering structure $\D: \R^2\rightrightarrows  \R^2$ defined by
$$
\D(y) := \begin{cases}
\D(y_1,y_1) & \mbox{ if }\; y_1 < y_2,\\
\D(y_2,y_2) & \mbox{ if }\; y_1 > y_2,\\
\conv\cone\left\{(1, \frac{1}{2}- \frac{a}{2|a|+1}), (\frac{1}{2}- \frac{a}{2|a|+1},1) \right\}  & \mbox{ if }\; y_1 = y_2= a.
\end{cases}
$$
It is easy to check that $\D(y) \subseteq \R^2_+$, that
$$
\forall a, b\in \R \; \mbox{ with } a\geq b,\;\D(a,a) \subseteq \D(b,b),
$$
and that condition (F3-a) is satisfied due to the convexity of the cones $\D(y)$ for every $y\in \R^2$. To check the fulfillment of (F3-c), fix two arbitrary elements $y,v \in \R^2$ such that $v\leq_N y$. It follows from the definition of $\D$ that
$$
y \in v + \D(v) \subseteq v + \R^2_+,
$$
which clearly implies that $y_1\geq v_1$ and $y_2 \geq v_2$. Denoting $\underline{y} := \min\{y_1,y_2\}$ and $\underline{v} := \min\{v_1,v_2\}$, we have $\underline{y} \geq \underline{v}$ and  $\D(\underline{y},\underline{y})\subseteq \D(\underline{v},\underline{v})$. Hence
$$
\D(y_1,y_2) = \D(\underline{y},\underline{y})\subseteq \D(\underline{v},\underline{v}) = \D(v_1,v_2).
$$
Since $\D(y_1,y_2)$ and $\D(v_1,v_2)$ are convex cones, we get $\D(v_1,v_2) + \D(y_1,y_2) \subseteq \D(v_1,v_2)$. Remembering that $y$ and $v$ were chosen arbitrarily in $\R^2$ allows us to conclude that $\D$ satisfies condition (F3-c), which completes our considerations in this example.}
\end{Example}

The last result of this section provides efficient conditions via the problem data that ensure the fulfillment of the limiting monotonicity assumption of Theorem~\ref{thm:EVP-dominated}. First we need to define the following notion of decreasing lower semicontinuity with respect to {\em nondomination relations}; cf.\ Definition~\ref{dec-lsc} in a different setting.

\begin{Definition}{\bf(decreasing lower semicontinuity with respect to nondomination)}\label{dec-lsc1} A vector-valued mapping $f: X\to Y$ is said to be {\sc $\leq_N$-decreasing lower semicontinuous with respect to nondomination relations} over a subset $A$ if for every forward-convergent sequence $\{x_n\} \subseteq A$ with a forward-limit $x_\ast$, we have that the decreasing monotonicity of $\{f(x_n)\}$ with respect to the nondomination relation $\leq_N$, i.e., $f(x_{n+1}) \leq_N f(x_n)$ for all $n\in \N$, implies that
$$
\forall n\in \N,\,f(x_\ast) \leq_N f(x_n).
$$
\end{Definition}

\begin{Proposition}[\bf sufficient conditions for limiting monotonicity in Theorem~\ref{thm:EVP-dominated}]\label{suff-limmon} Let $(X,q)$ be a forward-forward-complete quasimetric space, and let $f: X\to Y$ be a $\leq_N$-decreasing lower semicontinuous mapping with respect to the nondomination relation from Theorem~{\rm\ref{thm:EVP-dominated}}. Then the mapping $W$ defined by \eqref{def:W} satisfies the limiting monotonicity condition {\rm(F2)}.
\end{Proposition}
{\bf Proof}. Take an arbitrary generalized Picard sequence $\{x_n\}$ of the set-valued mapping $W$ satisfying $\sum_{n=0}^\infty q(x_n,x_{n+1})=L <\infty$. For each $\ve > 0$ there exists $N_\ve\in \N$ such that
$$
\sum_{n=N_\ve}^\infty q(x_n,x_{n+1}) > L -\ve/2,
$$
and thus for every $i, j \geq N_\ve$ with $j > i$ we have the estimates
$$
q(x_i,x_j) \leq \sum_{n=i}^{j-1} q(x_n,x_{n+1}) \leq \sum_{n=N_\ve}^\infty q(x_n,x_{n+1}) < \ve,
$$
which show that $\{x_n\}$ is a forward-Cauchy sequence. Since $X$ is forward-forward complete, this sequence forward-converges to some $x_\ast \in X$.

Remembering that $\{x_n\}$ is a generalized Picard sequence of $W$, we deduce from condition (F3-b) that
$f(x_{n+1}) \leq_N f(x_n)$ for all $n\in \N$. The imposed decreasing lower semicontinuity of $f$ gives us
$$
\forall n\in \N,\;f(x_{\ast}) \leq_N f(x_n).
$$
Using now (F3-c) yields the inclusions
\begin{eqnarray}\label{eqn:ss}
\forall n\in \N,\;\D(f(x_{n+1})) \subseteq \D(f(x_{n})) \; \mbox{ and } \; \D(f(x_\ast)) \subseteq \D(f(x_{n})).
\end{eqnarray}
Since  $\{x_n\}$ is a generalized Picard sequence of $W$, we get
\begin{eqnarray*}
f(x_{n+1}) \in  f(x_n) -\sqrt{\ve}q(x_n,x_{n+1})k - \D(f(x_{n+1})).
\end{eqnarray*}
Summing up the above relations from $n=i$ to $j$ while taking into account (F3-b) and the triangle inequality of the quasimetric ensure that
$$
f(x_{j+1}) \in f(i) - \sqrt{\ve}q(x_i,x_{j+1})k - \D(f(x_{j+1})).
$$
Adding the latter to the inclusion $f(x_\ast) \in f(x_{j+1}) - \D(f(x_{j+1}))$ yields
\begin{eqnarray*}
f(x_{\ast}) & \in&  f(x_i) - \sqrt{\ve}q(x_i,x_{j+1})k - \D(f(x_{\ast}))\\
& = & f(x_i) - \sqrt{q}(x_i,x_\ast)k + \sqrt{q}(x_{j+1}, x_\ast)k  \\
& & \hspace*{1in}- \sqrt{\ve}({q}(x_i,x_\ast) + {q}(x_{j+1}, x_\ast) - q(x_i,x_{j+1})k - \D(f(x_{\ast}))\\
& \st{(F3-b)}{\subseteq} & f(x_i) - \sqrt{q}(x_i,x_\ast)k + \sqrt{q}(x_{j+1}, x_\ast)k  - \D(f(x_{\ast})).
\end{eqnarray*}
Since the index $j$ was chosen arbitrarily and since $\D(f(x_{\ast}))$ is $k$-vectorial closed, the passage to limit above as $j\to\infty$ gives us the inclusion
$$
f(x_{\ast}) \in f(x_i) - \sqrt{q}(x_i,x_\ast)k - \D(f(x_{\ast})),
$$
i.e., $x_\ast \in W(x_i)$. This verifies (F2) by taking into account that $i$ was also chosen arbitrarily. $\h$\vspace*{-0.1in}

%%%%%%%%%%%%%%%%%%%%%%%%%%%%%%%%%%%%%%%%%%%%%%%%%%%%%%%%%%%%%%%%%%%%%%%%%%

\section{Variational Rationality in Behavioral Sciences}\label{behav}

This section is devoted to developing a {\em variational rationality approach} to human dynamics in the vein of \cite{s09,s10,s16,s19a,s19b}, with taking now into account the new results on variational principles in vector optimization with variable domination structures that were established above.

After presenting the basic concepts of the variational rationality modeling of human dynamics, we mainly concentrate on the following issues:\\[1ex]
$\bullet$ Introducing {\em generalized efficiency} and {\em domination structures} of the type formalized in Definition \ref{def:binary}, while being adjusted to the variational rationality approach to human dynamics. These notions extend those from  \cite{yu74} to the settings where the {\em resistance to move matters}.\\[1ex]
$\bullet$ Applying the obtained variational principles in vector optimization with variable domination structures to establish the existence of {\em ex ante} (before moving) and {\em ex post} (after moving) variational traps with showing that {\em possible regrets can matter much}.\vspace*{0.05in}

To highlight the major topics of the presentation, we split this section into several subsections.

\subsection{Pareto and Yu Efficiency and Domination Structures in Behavioral Models}

In this subsection we discuss some basic notions in the modeling of human dynamics and conventional approaches to behavioral models based on efficiency and domination structures in the classical sense of Pareto and more recent ones introduced by Yu.\\[1ex]
{\bf 5.1.1: Pareto Efficiency and Nondomination Binary Relations in Behavioral Models.} In the finite-dimensional space $Y=\mathbb{R}^m$, consider a list of different {\em pains} $J =\left\{1,\ldots,m\right\}$ and a list of vector {\em amounts of pains} $v=(v^{1},\ldots,v^{m})\in\R_{+}^{m} \subseteq Y$ associated with each pain $j\in J$. The amount of pain $v^{j}\in\R_{+}$ represents a quantity of a {\em ``to be decreased" payoff}, e.g., some degree of unsatisfaction, loss, cost, lack of given things (size of needs), etc. In this context where the agent wants less of each pain and tries to minimize the amounts of pain $v^{j}$ as $j\in J$, the space $Y=\R^{m}$ is the {\em space of amounts} associated with pains and pleasures. The first problem of the agent is to compare the lists of different amounts of pains $v = (v^{1},\ldots,v^{m})\in\R_{+}^{m}$ and $y=(y^{1},\ldots,y^{m})\in\R_{+}^{m}$. As usual, a vector of pains $v$ is {\em Pareto smaller} (resp.\ {\em Pareto larger}) than another vector of pains $y$ if and only if we have $v^{j}\leq y^{j}$ (resp.\ $v^{j}\geq y^{j}$) for all $j\in J$.

Having in mind the above descriptions of ``Pareto smaller" and ``Pareto larger" vectors, we need now to clarify the related meanings of the expressions that a given vector of pains is
``better than" or is ``worse than" another vector of pains in all their aspects; cf.\ Definition~\ref{def:Pareto}.

\begin{description}
\item[(Pareto-i)] The vector of pains $v$ is {\em Pareto better} than the vector of pains $y$ (in all their aspects) if $v=y-D^{\ast }$, where $D^*:=\R^m_+$. In this setting, $v=y-d$ with $d\in D^{\ast}$ means that each amount of pain $v^{j}$ of the list $v$ is smaller than or equal to each amount of pains $y^{j}$ of the list $y$, i.e., $v^{j}\leq y^{j}$ for all $j\in J$.
This tells us that ``less of each pain is better."

\item[(Pareto-ii)] The vector of pains $y$ is {\em Pareto worse} than the vector of pains $v$ (in all their aspects) if $y\in v + D^{\ast }$, where $D^{\ast}:=\R_{+}^{m}$. In this
case we have $y=v+d$ with $d\in D^{\ast}$ meaning that each amount of pain $y^{j}$ of the list $y$ is higher or equal to each amount of pains $v^{j}$ of the list $v$, i.e., $y^{j}\geq v^{j}$ for all $j\in J$. This tells us that ``more of each pain is worse."
\end{description}

Note that in the above case of $D^{\ast}:=\R_{+}^{m}$ we have the equivalencies (cf.\  Section~\ref{sec:domin})
$$
y-v \in D^{\ast } \Longleftrightarrow v\leq_{D^{\ast}} y\Longleftrightarrow v-y\leq _{D^{\ast }}0.
$$
Thus defining the ``better than" sets $\mathfrak{B}(y)=y-D^{\ast }$ and the ``worse than" sets $\mathfrak{W}(y)=y+D^{\ast }$, we get that $v$ is Pareto better (resp.\ worse)
than $y$ {\em if and only if} $v\in\mathfrak{B}(y)$ (resp.\ $y\in \mathfrak{W}(y)$).\\[1ex]
{\bf 5.1.2: Yu Efficiency and Nondomination Binary Relations} (cf.\ Definitions~\ref{def:binary} and \ref{def:sol}). In \cite{yu73dom,yu74}, Yu generalized the above concepts of
Pareto efficiency and nondomination by considering the following {\em conic domination structures}:
\begin{description}
\item[(Yu-i)] An arbitrary fixed convex cone $D\subseteq Y$ instead of the Pareto constant cone $D^{\ast}=\mathbb{R}^m$.
\item[(Yu-ii)] Variable conic structures $\mathcal{D}(y)$ for all $y\in Y$.
\end{description}

To discuss these concepts, consider first the ``worse" and ``better" relations with respect to constant cones $D \subseteq Y$. The main emphasis here is that the same
amount of two different pains can be more important or less important for an individual. In this setting, the agent may accept to trade off a lower amount of a more important pain $1$ to a higher amount of a less important pain $2$. For simplicity, take two pains $1$ and $2$ with the corresponding amounts of these pains $v=(v^{1},v^{2})\in\R_{+}^{2}$ and $y=(y^{1},y^{2})\in\R_{+}^{2}$. The new meaning that Yu gave to the ``better than" relation is: $v$ is better than $y$ if and only if $v=y-D$, when the cone $D$ may be larger than the Pareto cone $D^{\ast}=\R_{+}^{m}$. In this simple situation where the agent compares the new short list of pains $v$ with the old short list of pains $y$, he/she prefers the new list of pains $v$ to the old one $y$ if accepting to trade off the lower amount $v^{1}\leq y^{1}$ of the most important pain $1$ against a higher amount $\ v^{2}\geq y^{2}$ of the less important pain $2$. This means that, moving from the list $y$ to the list $v\in y-D$, the agent trades off the diminution $v^{1}-y^{1}\leq 0$ of the most important pain against the augmentation $v^{2}-y^{2}\geq 0$ of the less important pain. However, this augmentation should not be too important relative to the diminution. Given the diminution  $y^{1}-v^{1}=d^{1}\geq 0$ of pain~1, the larger augmentation $v^{2}-y^{2}=-d^{2}\geq 0$ of pain~2, which the agent tolerates, can be
$$
0\leq v^{2}-y^{2}=-d^{2}\leq \alpha ^{1}(y^{1}-v^{1})=\alpha ^{1}d^{1}
$$
giving us $d^{2}+\alpha ^{1}d^{1}\geq 0$. Then, $\alpha^{1}>0$ is the larger augmentation of the less important pain~2 that the agent would accept relative
to the diminution of one unit of the more important pain~1, i.e.,
$$
0\leq v^{2}-y^{2}\leq \alpha^{1}\;\mbox{ if }\;y^{1}-v^{1}=1.
$$
These considerations lead us to the following construction of the cone $D$ of {\em acceptable augmentations} of less important pains relative to given
diminutions of more important pains:
$$
D =\big\{d=(d^{1},d^{2})\in\R^{2}\;\big|\;d^{1}\geq 0,\;\alpha ^{1}d^{1}+d^{2}\geq 0\big\}.
$$
We also refer the reader to \cite{hwh10} for related discussions showing how tradeoffs modelize the relative importance of criteria in the case where $D\supseteq\R_{+}^{m}$ is not a too obtuse cone.

The same meaning can be given to the relation ``to be worse than" with the same dominated cone $D$. In this case we say that $y$ is worse than $v$ if $y\in v+D$, i.e.,
if $y=v+d$ as $d\in D$. This tells us that the agent finds $y$ worse than $v$ if a given augmentation of the more important pain~1 is not compensated by a large enough diminution of the less important pain~2.

\subsection{Variational Rationality with Variable Efficiency and Domination Structures}

In this subsection we describe the {\em variational rationality {\rm(VR)} approach} to human dynamics, which benefits from the variational principles developed in Section~\ref{new}. In fact, one of the major motivations for developing our new research on variational principles in vector optimization problems with variable domination structures came from the needs of the VR approach described below.\\[1ex]
{\bf 5.2.1: Ex Ante and Ex Post Visions of Moves.} According to Definition~\ref{def:binary}, a variable domination and efficiency structures gives us, for each position $v\in Y$ of a space of pains (positions) $Y$, a set of possible worse (dominated) positions $v+\mathcal{D}(v)\subseteq Y$ and a set of better (preferred) positions $v-\mathcal{D}(v)\subseteq Y$. The variable cone $\mathcal{D}(v)\subseteq Y$ represents a set of pains, which can be added to the vector of pains $v$ to make worse the new vector of pains $v+d$ with $d\in\mathcal{D}(v)$. On the other hand, the variable cone $\mathcal{P}(v) = -\mathcal{D}(v)\subseteq Y$ defines a set of pains that can be dropped from the vector of pains $v$ while making better the vector of pains $v-d$ with $d\in \mathcal{D}(v)$.\vspace*{0.03in}

The {\em variable nondomination and efficiency binary relations} under consideration are given by:
{\begin{itemize}
\item ``$y$ is worse than $v$" if and only if  $y \in v+\mathcal{D}(v)$; this is the nondomination binary relation taken from Definition~\ref{def:binary}(i).
\end{itemize}
\begin{itemize}
\item
``$v$ is better than $y$" if and only if $v \in y-\mathcal{D}(y)$; this is the efficiency binary relation taken from Definition~\ref{def:binary}(ii).
\end{itemize}

Note that these variable relations are not generally equivalent, while they becomes equivalent in the case of constant structures $\mathcal{D}(y)=\mathcal{D}(v)=D$ for all $v,y\in Y$.\vspace*{0.03in}

The VR approach to human dynamics focuses the major attention on a short list of main concepts for modeling human behaviors: activities, payoffs (utilities and disutilities as satisfactions and unsatisfactions to move), moves, costs to move, advantages/disadvantages to move, inconveniences to move, motivation and resistance to move, worthwhile balances, worthwhile moves, aspiration points, desires, and stationary or variational traps. This approach is well adapted to: (a) give a new interpretation of the Yu's approach in the context of variable cones when there are no resistance to
move (i.e., change rather than stay), and (b) generalize the nondomination and efficiency binary relations in the vein of Definition~\ref{def:binary} when resistance to move matters. Without resistance to move, the VR approach provides the following:
\begin{description}
\item[(i)] Starts by focusing the attention on moves in a space of positions.
\item[(ii)] Makes a distinction between disadvantageous (utility deteriorating) moves $d\in\mathcal{D}(v)$ and advantageous (utility improving) moves $-d$ with $d\in \mathcal{D}(v)$ in the space of positions.
\item[(iii)] Makes an essential distinction between an ex ante perception to move and an ex post perception to move. In this context, $\mathcal{D}(v)$ is a set of disadvantageous moves, while $\mathcal{P}(v)\ = -\mathcal{D}(v)$ is a set of advantageous moves starting from the initial position $v$ in the payoff space $Y$.
\end{description}

We refer the reader to \cite{bms15,bms15b} for the first attempts to investigate adaptive aspects of the variational rationality approach when resistance to move matters. Now we can do more.\\[1.3ex]
{\bf 5.2.2: Should I stay, should I move?} The corresponding logic of the efficiency and domination structures become very clear in the context of the variational rationality approach to human
dynamics. We start here with the VR discussions concerning the space of ``to be decreased" payoffs $Y=\R^{m}$. Let $y\in Y$ and $v\in Y$ be the amounts of pains that the agent endorses in the previous
and the current periods. Then, within the current period, a simplified definition of a move that is well adapted to the present paper starts with ``having suffered of the amounts of pains $y\in Y$ in the previous period" and ends with ``suffering of the amount of pains $v$" in the current period. This move is $(y,v)\in Y\times Y$. It is a change if $v\neq y$ and a stay if $v=y$. Note that the move $(y,v)$ in the payoff space corresponds to some move $(x,u)\in X\times{X}$ in the activity space $X$. With the two given bundles of activities $x\in X$ and $u\in X$, we have the amount of pains $y = f(x)\in Y$ in the previous period and the amounts of pains $v = f(u)\in Y$ in the current period. The most basic question driving the VR approach is the following: ``should I stay or should I move?" That is, at the beginning of the current period (ex
ante, i.e., before moving) the main alternative is:
\begin{itemize}
\item[\bf(a)] Either {\em to stay}, i.e., doing the same bundle of activities $x$ in the current period as before. In this case, the agent would suffer from the same amounts of
pains $y=f(x)$ as before.

\item[\bf(b)] Or {\em to change}, i.e., doing a different bundle of activities $u\neq x$ in the current period as before. In this case, the agent will suffer from new amounts of pains $v=f(u)$ in the current period.
\end{itemize}

Let us discuss the aforementioned major alternative from both viewpoints of the {\em efficiency} and {\em nondomination} binary relations introduced in Definition~\ref{def:binary}.\\[1ex]
{\sc Efficiency binary relation} (Definition~\ref{def:binary} (ii))}: {\sc should I change?} {\bf Yes, if ex ante $v$ is better than $y$.} In this case the advantages to move from $y$ to $v$
(change rather than stay) in the payoff space is $\mathbb{A}(y,v):=y-v = f(x)-f(u) = A(x,u) \in Y$. Consider an ex ante perception of a move. In this setting, the agent prefers to change before
moving from $x$ to $u$, rather than to stay at $x$, if the new amount of pains $v = f(u)$ is lower than the old one $y = f(x)$. This means that, ex ante, a given diminution of the most important pains compensates a not too large augmentation of the less important pains. The latter is equivalent to saying that $\mathbb{A}(y,v) = y-v \in\mathcal{D}(y) \Longleftrightarrow A(x,u) = f(x)-f(u) \in \mathcal{D}(f(x))$, which means that there are ex ante advantages to move from $y$ to $v$, i.e., from $x$ to $u$. Appealing to Definition~\ref{def:binary}(ii)), this can be written as $v \in y- \mathcal{D}(y)$.\\[1ex]
{\sc Nondomination binary relation} (Definition~\ref{def:binary}(i)): {\sc should I regret to have changed?} {\bf No, if ex post $y$  is worse than $v$.} Indeed, consider, the agent's ex post perception of the same move $(y,v)$. In this new setting, the agent would prefer to change from $y$ to $v$ after moving, i.e., to go from $x$ to $u$ rather than to stay at $y$ provided that the new amount of pains $v = f(u)$ is perceived ex post as lower than the old amount of pains $y = f(x)$. This means that after moving a given diminution of the most important pains compensates a not too large augmentation of the less important pains. The latter is equivalent to saying that $\mathbb{A}(y,v): = y-v \in \mathcal{D}(v)\Longleftrightarrow A(x,u) = f(x)-f(u) \in \mathcal{D}(f(u))$, which tells us that there are ex post advantages to move from $y$ to $v$, i.e., from $x$ to $u$. Coming back to Definition~\ref{def:binary}(i), this can be written $y \in v + \mathcal{D}(v)$ meaning that the agent does not regret   ex post to move from $y$ to $v$. We refer the reader to \cite{s19b} for more discussions of the possible origins of regrets in the variational rationality context, where ex post regrets come from wrong ex ante evaluations of utility and costs of different moves.\\[1.ex]
{\bf 5.2.3: When the resistance to move matters much}. First let us offer an appropriate extension of variable domination and efficiency structures in the variational rationality approach when the {\em resistance to move matters}. With respect to the binary relations in Definition~\ref{def:binary}, the VR approach compares advantages to move to inconveniences/resistance to move and defines ex ante and ex post {\em worthwhile moves} that generalize ex ante and ex post {\em advantageous moves}, respectively.

\textbf{Inconveniences to move.} When the resistance to move (change rather than stay) {\em matters}, the resistance to move generates the {\em inconvenience to move rather than to stay} defined by $\mathbb{I(}y,v):= \mathbb{C}(y,v) - \mathbb{C}(y,y) = C(x,u)-C(x,x) = I(x,u)$. In this formula, the amount $\mathbb{C}(y,v) = C(x,u) \in\R_{+}^{m} \subseteq Y$ represents vectorial costs to move from ``having done the bundle of activities $x$ in the previous period" to ``being able to do and do the bundle of activities $u$ in the current period". The amount $\mathbb{C} (y,y) = C(x,x) \in\R_{+}^{m}\subseteq Y$ defines vectorial costs to stay at $x$. Observe that costs to move are not symmetric, i.e., $\mathbb{C}(v,y)\neq \mathbb{C}(y,v)$. This requires to use in modeling the framework of {\em quasimetric spaces}, which has been done in the variational theory developed in Section~\ref{new}. In our VR behavioral model we consider a specific vectorial case, where $I(x,u)=\sqrt{\varepsilon }q(x,u)k$  with $q(x,u)$ being a given quasidistance; see \cite{bms15,bms15b} for more discussions of such issues.

\textbf{Advantageous moves}. When resistance to move {\em does not matter}, we define in our terminology (following an implicit construction of \cite{s19b}) an {\em advantageous move} $(y,v)$ from the viewpoint of $r\in \left\{ y,v\right\}$ by $\mathbb{A}(y,v):=y-v \in \mathcal{D}(r) \Longleftrightarrow A(x,u) = f(x)-f(u) \in
\mathcal{D}(r)$.

\textbf{Worthwhile moves}. Now we are ready to define, based on the binary relations from Definition~\ref{def:binary}, the new notion of worthwhile moves when the {\em resistance to move matters}. Namely, the {\em worthwhile move} $(y,v)$ from the viewpoint of $r\in\left\{y=f(x),v=f(u),y_{0} = f(x_{0})\right\}$ is given by $\mathbb{B}_{\xi }(y,v): = \mathbb{A}(y,v)-\xi
\mathbb{I(}y,v) \in \mathcal{D}(r)$. In our specific context of efficient and domination structures, the {\em worthwhile balance} is defined by
\begin{equation*}
B_{\xi }(y,v): = \mathbb{A}(y,v) - \xi \mathbb{I}(y,v) = y-v-\xi \mathbb{I}(y,v), \;\mbox{ with } \; \xi >0.
\end{equation*}
In \cite{s09,s10,s16}, the reader can find some discussions on motivation of the resistance to move in other behavioral science settings. Note that the concepts from Definition~\ref{def:binary} correspond to the case where $\xi:=\sqrt{\varepsilon} =0$. In what follows we consider even a more specific balance situation with $\mathbb{I}(y,v) = \sqrt{\varepsilon}q(x,u)k$, i.e.,
$$
B_{\xi}(y,v) = f(x) - f(u) - \sqrt{\varepsilon}q(x,u)k,
$$
where $\xi = \sqrt{\varepsilon }\geq 0$. Then, a move $(y,v)$, which starts from the position $y$ and goes to the position $v$, is {\em worthwhile} in the following senses:\\[1ex]
$\bullet$ {\em Ex ante} if we have $B_{\xi }(y,v)\in \mathcal{D}(y)$ before moving, while choosing the viewpoint of the starting position $r=y=f(x)$ and the viewpoint of $x$ (efficiency binary relation as in Definition~\ref{def:binary}(ii)).\\[1ex]
$\bullet$ {\em Ex post} if we have $B_{\xi }(y,v)\in \mathcal{D}(v)$ after moving, while choosing the viewpoint of the final position $r=v=f(u)$ and the viewpoint of $u$ (nondomination binary relation as in Definition~\ref{def:binary}(i)).\vspace*{0.05in}

Observe finally that the move $(y_{0} = f(x_{0}),v=f(u))$ is {\em ex ante worthwhile} if
$$
B_{\xi }(y_{0},v)=\mathbb{A}(y_{0},v)-\xi \mathbb{I(}y_{0},v)=f(x_{0})-f(u)-\sqrt{\varepsilon }q(x_{0},u)k\in \mathcal{D}(v_{0}),
$$
while choosing the initial viewpoint of $r=y_{0} = f(x_{0})$ and from viewpoint of $x_{0}$.

\subsection{Existence of Ex Ante and Ex Post Variational Traps}

The above discussions show that the behavioral model of human dynamics, which is described in terms of the variational rationality approach, can be enclosed into the variational framework of vector optimization with variable domination structures in quasimetric spaces. Then, the new variational principles of Section~\ref{new} obtained in this general framework leads us to behavioral conclusions that can be interpreted as the {\em existence of ex ante} and {\em ex post variational traps}. The results presented below are direct consequences of the obtained variational principles, which are derived in Theorems~\ref{thm:Pareto} and \ref{thm:EVP-dominated}. Note that, besides the statements of the these theorems, their very {\em proofs} based on constructive generalized Picard sequences provide efficient {\em dynamic procedures} to approach such traps, not only to establish their existence.\vspace*{0.05in}

To this end, we define in the framework of the VR approach the concept of variational traps as follows. A given position $x_{\ast}$ is a {\em variational trap} if this position is worthwhile to reach, but not worthwhile to leave. Using the definition of a worthwhile balance formulated in this section and the results obtained in Theorems~\ref{thm:Pareto} and \ref{thm:EVP-dominated} allows us to arrive at the following conclusions:\\[1ex]
{\bf $\bullet$ Ex ante variational traps}. It follows from the results of Theorem~\ref{thm:Pareto} with $\xi:=\sqrt{\ve}>0$ that we have the conditions:
\begin{itemize}
\item[\bf(i)] $B_{\xi}(y_{0}=f(x_{0}),v_{\ast }=f(x_{\ast }))=f(x_{0})-f(x_{\ast })-\sqrt{\varepsilon }q(x_{0},x_{\ast })k\in \mathcal{D}( f(x_{0}))$.
\item[\bf(ii)] $B_{\xi }(v_{\ast }=f(x_{\ast }),y=f(x))=f(x_{\ast })-f(x)-\sqrt{\varepsilon }q(x_{\ast },x)k\in \mathcal{D}(f(x_{0}))$.
\end{itemize}

Note that condition (i) means that it is worthwhile to move from $x_{0}$ to $x_{\ast}$, while condition (ii) tells us that it is not worthwhile to move away from $x_{\ast}$. The point of view that determines preferences is the initial position $(x_{0},r=f(x_{0}))$. This defines an ex ante variational trap as an efficiency binary relation from Definition~\ref{def:binary}(ii)), which gives us the
ex ante  motivation to move from $x_{0}$ to $x_{\ast }$ and then to stay at $x_{\ast}$.\\[1ex]
{\bf $\bullet$ Ex post variational traps}. It follows from the results of Theorem~\ref{thm:EVP-dominated} with $\xi:=\sqrt{\ve}>0$ that we have the conditions:
\begin{itemize}
\item[\bf(i)] $B_{\xi }(y_{0}=f(x_{0}),v_{\ast }=f(x_{\ast }))=f(x_{0})-f(x_{\ast })-\sqrt{\varepsilon }q(x_{0},x_{\ast })k\in \mathcal{D}( f(x_{\ast}))$.
\item[\bf(ii)] $B_{\xi }(v_{\ast }=f(x_{\ast }),y=f(x))=f(x_{\ast })-f(x)-\sqrt{\varepsilon }q(x_{\ast },x)k\in \mathcal{D}(f(x))$.
\end{itemize}
As seen, condition (i) tells us that it is worthwhile to move from $x_{0}$ to $x_{\ast}$, while condition (ii) means that it is not worthwhile to move away from $x_{\ast}$. The point of view that determines preferences in this case for condition (i) is the final position $(x_{\ast},r=f(x_{\ast }))$. On the other hand, for condition (ii) it is the position $x$ with $r=f(x)$ for each $x$ away from $x_{\ast }$. This defines an ex post variational trap corresponding to the nondomination binary relation from Definition~\ref{def:binary}(i), which excludes ex post regrets to move from $x_{0}$ to $x_{\ast}$ and then to stay at $x_{\ast}$.

Thus Theorems~\ref{thm:Pareto} and \ref{thm:EVP-dominated} provide efficient conditions ensuring the existence of ex ante and ex post variational traps in the variational rationality model of human dynamics.\vspace*{0.05in}

{\bf Acknowledgements}. The research of the first author was initially conducted during his stay at the Vietnam Institute for Advanced Study in Mathematics (VIASM), Hanoi, Vietnam. He would like to
thank the institute for hospitality and support.

\end{document}